\documentclass[preprint]{elsarticle}

\usepackage{lineno,hyperref}
\usepackage{indentfirst,latexsym,bm}
\usepackage{mathrsfs}
\usepackage{CJK,fancybox,fancyhdr,color,graphicx,amsmath,amsthm,amssymb,enumerate,amsfonts}
\usepackage{anysize}
\usepackage[misc]{ifsym}
\usepackage[T1]{fontenc}
\usepackage{subfigure}
\usepackage{algorithmicx,algorithm, algpseudocode}
\usepackage{xcolor}
\usepackage{colortbl,booktabs}

\newtheorem{thm}{Theorem}[section]
\newtheorem{cor}[thm]{Corollary}
\newtheorem{lem}[thm]{Lemma}

\newtheorem{defn}[thm]{Definition}
\newtheorem{rem}[thm]{Remark}

\modulolinenumbers[5]

\journal{Journal of \LaTeX\ Templates}

\makeatletter

\newenvironment{breakablealgorithm}
{
\begin{center}
\refstepcounter{algorithm}
\hrule height.8pt depth0pt \kern2pt
\renewcommand{\caption}[2][\relax]{
{\raggedright\textbf{\ALG@name~\thealgorithm} ##2\par}%
\ifx\relax##1\relax 
\addcontentsline{loa}{algorithm}{\protect\numberline{\thealgorithm}##2}%
\else 
\addcontentsline{loa}{algorithm}{\protect\numberline{\thealgorithm}##1}%
\fi
\kern2pt\hrule\kern2pt
}
}{
\kern2pt\hrule\relax
\end{center}
}

\makeatother

\begin{document}

\begin{frontmatter}

\title{The standard forms and convergence theory of the Kaczmarz-Tanabe type methods for solving linear systems}

\author{Chuan-gang Kang}
\address{School of Mathematical Sciences, Tiangong University, Tianjin 300387, People's Republic of China}
\ead{ckangtj@tiangong.edu.cn}



\begin{abstract}
    In this paper, we consider the standard forms of two kinds of Kaczmarz-Tanabe type methods, one is derived from the Kaczmarz method and the other is derived from the symmetric Kaczmarz method. As a famous image reconstruction method in computerized tomography, the Kaczmarz method is simple and easy to implement, but its convergence speed is slow, so is the symmetric Kaczmarz method. When the standard forms of the Kaczmarz-Tanabe type methods are obtained, their iteration matrices can be used continuously in the subsequent iterations. Moreover, the iteration matrices can be stored in the image reconstruction devices, which enables the Kaczmarz method and the symmetric Kaczmarz method to be used like the simultaneous iterative reconstructive techniques (SIRT). Meanwhile, theoretical analysis shows that the convergence rate of the symmetric Kaczmarz-Tanabe method is better than that of the Kaczmarz-Tanabe method but is slightly worse than that of two-step Kaczmarz-Tanabe method, which is verified numerically. Numerical experiments also show that the convergence rates of the Kaczmarz-Tanabe method and the symmetric Kaczmarz-Tanabe method are better than those of the SIRT methods.
\end{abstract}


\end{frontmatter}


\section{Introduction}\label{section.introduction}

In medical imaging tomography (see, i.e., \cite{Ledley1988,Herman2010Fundamentals,Natterer2001}), people are often asked to solve the following linear system of equations, i.e.,
\begin{align}\label{linear.system}
  Ax=b,
\end{align}
where $A\in\mathbb{R}^{m\times n}, b\in\mathbb{R}^m$ are also called projection matrix and measurement vector, respectively. We suppose that \eqref{linear.system} is consistent and $x^*$ is a true solution. If $A$ is not full column rank, $x^{\dagger}=A^\dagger b$ is used to denote the minimum norm least-squares solution \cite{Engl96,Wang2021Gauss} of \eqref{linear.system}, where $A^\dagger$ denotes the pseudo-inverse of $A$. The linear system \eqref{linear.system} can be generated by discretizing the Radon transform
\begin{align*}
  p=\int_L f(x,y)ds,
\end{align*}
where, $L(\rho,\theta)=\{(x,y): x\cos\theta+y\sin\theta=\rho\}$ is the path of integration, $f(x,y)$ is the relative attenuation of the object to ray at point $(x,y)$ on the line $L$ and $ds=\sqrt{\rho^2+\rho'(\theta)^2}d\theta$; let $\phi$ denote the angle between the normal direction of $L$ and the polar axis on a given complex plane, so $\theta\in (\phi-\pi/2,\phi+\pi/2)$, (see, e.g., \cite{Herman2010Fundamentals,Radon1917,Herman1995,Gordon1970}).

The Kaczmarz method proposed by the Polish mathematician Kaczmarz\cite{Kaczmarz1937} is one of the most popular iterative methods to solve \eqref{linear.system} in computerized tomography. Let $A=(a_1,a_2,\ldots,a_m)^T$, then the Kaczmarz's iteration reads
\begin{align}\label{Kaczmarz.iteration}
  x_k=x_{k-1}+\frac{b_i-\langle a_i,x_{k-1}\rangle}{\|a_i\|_2^2}a_i,\quad k=1,2,\ldots,
\end{align}
where $i=\bmod(k-1,m)+1$, $\langle x,y\rangle=x^Ty$ and $\|x\|_2=\sqrt{\langle x,x\rangle}$ denote the inner product of $x,y$ and the $2$-norm of $x$ in $\mathbb{R}^n$, respectively.

The symmetric Kaczmarz method can be described as
\begin{align}\label{symmetric.Kaczmarz.orig}
  x_k=x_{k-1}+\frac{b_i-\langle a_i,x_{k-1}\rangle}{\|a_i\|_2^2}a_i,\quad k=1,2,\ldots,
\end{align}
where
\begin{align}\label{symmetric.Kaczmarz.con}
  &i=\left \{
    \begin{array}{ll}
       \text{mod}(k,2m-2),    & 1\le\text{mod}(k,2m-2)\le m,\\
       2m-\text{mod}(k,2m-2), & m<\text{mod}(k,2m-2)\le 2m-3,\\
       2, & \text{mod}(k,2m-2)=0.\\
    \end{array}
  \right.
\end{align}
Compared with the popular expression of the symmetric Kaczmarz method (see, i.e., \cite{Elfving1979,Huang1993The}), the iterative scheme \eqref{symmetric.Kaczmarz.orig} is more consistent in form with Kaczmarz's iteration.

The Kaczmarz method has many advantages, such as good convergence, ease to implement and so on, and has been used to solve the phase problem \cite{Wei2015}. However, the convergence speed of the Kaczmarz method sometimes becomes very slow, especially when the successive hyperplanes meet at a very small angle. In order to keep the advantages of the Kaczmarz method and overcome its disadvantages, many scholars consider the subsequence $\{y_k\}$ of sequence $\{x_k\}$, where $y_k=x_{k\cdot m}$. Kang \cite{Kang2021Convergence} gives the following iterative scheme of Kaczmarz's subsequence $\{y_k\}$, i.e.,
\begin{align}\label{Kaczmarz.Tanabe.iteration.orig}
  y_{k+1}=(I-A_{\mathcal{S}}^TMA)y_k+A_{\mathcal{S}}^TMb,\qquad k=0,1,2,\ldots,
\end{align}
where $I$ denotes the identity matrix of whatever size appropriate to the context, and
\begin{align}
  &P_i=I-\frac{a_ia_i^T}{\|a_i\|_2^2},\quad i=1,2,\ldots,m,\\
  &Q_m=I, Q_j=P_mP_{m-1}\ldots P_{j+1}, \quad j=1,2,\ldots,m-1, \label{definition.Qj}\\
  &Q=P_mP_{m-1}\cdots P_1,\label{definition.Q}\\
  &A_\mathcal{S}=(Q_1a_1,Q_2a_2,\ldots,Q_ma_m)^T,\label{definition.AS}\\
  &M=\text{diag}(1/\|a_1\|_2^2,1/\|a_2\|_2^2,\ldots,1/\|a_m\|_2^2)\label{definition.M}.
\end{align}
and the subsequent iteration \eqref{Kaczmarz.Tanabe.iteration.orig} was named the Kaczmarz-Tanabe's iteration by Popa \cite{Popa2018}.

Compared with Kaczmarz's iteration, Kaczmarz-Tanabe's iteration has good approximate stability (i.e., iterative error does not fluctuate as violently as Kaczmarz's iteration (see \cite{Kang2021Convergence}), which may provide convenience for people to study the regularization theory of the Kaczmarz method). In fact, compared with the traditional iterative scheme of Kaczmarz-Tanabe method (see \cite{Popa2018, Tanabe1971}), there are many improvements in the expression of \eqref{Kaczmarz.Tanabe.iteration.orig} . However, $A_s$ is the compound of $Q_i$ and $a_i$, which brings many obstacles for further research, especially the regularization theory, etc. In this paper, we mainly consider the standard form of \eqref{Kaczmarz.Tanabe.iteration.orig}, and the corresponding iteration matrix can be calculated by blocking and parallelization techniques.

Assume that $\text{rank}(Q)=p$, and $\sigma_1,\sigma_2,\ldots, \sigma_p$ are $p$ non-zero singular values of $Q$. Kang gave the following convergence result (see \cite[Theorem 2.10 \& Corollary 2.11]{Kang2021Convergence}).

\begin{thm}\label{thm.convergence.rate.Kaczmarz.Tanabe}\cite{Kang2021Convergence} For any matrix $A$ without zero row, let $\{y_k,k\ge 0\}$ be the sequence of vectors generated by \eqref{Kaczmarz.Tanabe.iteration.orig} and $e_k:=y_k-x^\dagger-P_{N(A)}y_0$, then
\begin{align*}
  \|e_{k+1}\|_2\le \max_{0<\sigma_i<1}\sigma_i\|e_k\|_2,\quad \|e_{k+1}\|_2\le \max_{0<\sigma_i<1}\sigma_i^{k+1}\|e_0\|_2
\end{align*}
holds, where $\sigma_i$ is the singular value of $Q$.
\end{thm}

We next consider the Kaczmarz-Tanabe method and hope to get a matrix-vector form similar to the SIRT methods. For ease of reference, we list several typical representations of the SIRT methods (see, i.e., \cite{Hansen2018AIRtools,Elfving2010}) and the general iteration reads
\begin{align}\label{the.SIRT.method}
  x_{k+1}=x_k+\lambda_kTA^TM(b-Ax_k),
\end{align}
where $\lambda_k$ is the relaxation parameter. For each $j=1,2,\ldots,n$, we denote by nz$_j$ the number of nonzero elements in the $j$-th column of A, and $S=\text{diag}(\text{nz}_1,\ldots,\text{nz}_n)$. For convenience of description, we also denote the sum of the $i$-th row of $A$ by $s_{r_i}$ and the sum of the $j$-th column of $A$ by $s_{c_j}$. Let $\|x\|_S=\sqrt{x^TSx}$ denote a weighted Euclidean norm. When $\lambda_k\equiv 1$, the following methods will be obtained by taking given $T$ and $M$ pairs.
\begin{itemize}
  \item Landweber\cite{Landweber51}:\quad $T=I, M=I$;
  \item Cimmino\cite{Cimmino1938}: \quad $T=I, M=D=\frac{1}{m}\text{diag}(\frac{1}{\|a_1\|_2^2},\ldots,\frac{1}{\|a_m\|_2^2})$;
  \item CAV\cite{Censor2001Component}: \qquad\quad $T=I, M=D_S=\text{diag}(\frac{1}{\|a_1\|_S^2},\ldots,\frac{1}{\|a_m\|_S^2})$;
  \item DROP\cite{Censor2008On}: \qquad $T=S^{-1}, M=mD$;
  \item SART\cite{Jiang2003Convergence,Andersen1984Simultaneous,Wan2011Three}: \quad  $T=\text{diag}(s_{c_1},\ldots,s_{c_n})^{-1}, M=\text{diag}(s_{r_1},\ldots,s_{r_m})^{-1}$.
\end{itemize}

The rest of the work is organized as follows. In Section \ref{section.further.study.Kaczmarz.Tanabe method}, we consider the standard form (i.e., matrix-vector form) of the Kaczmarz-Tanabe method, and introduce some concepts related with the sequential projection. In Section \ref{section.the.symmetric.Kacmarz.Tanbe}, we consider the matrix-vector form of the symmetric Kaczmarz method and analyze its convergence rate. In section \ref{section.algorithm}, we give the algorithm flows to calculate $C$ appearing in \eqref{Kaczmarz.Tanabe.iteration.standard} and $\bar{C}$ appearing in \eqref{symmetric.Kaczmarz.Tanabe.iteration.specific} respectively. In Section \ref{section.numerical.tests}, we compare the computational efficiency of the Kaczmarz-Tanabe method, symmetric Kaczmarz-Tanabe method, SIRT methods and CGMN method \cite{Elfving1979} by numerical experiments.

\section{The standard form of the Kaczmarz-Tanabe method and its convergence}
\label{section.further.study.Kaczmarz.Tanabe method}

Compared with \eqref{Kaczmarz.iteration}, the Kaczmarz-Tanabe iteration \eqref{Kaczmarz.Tanabe.iteration.orig} has made great change in form because it gets rid of the constraint of projection row by row according to the system of equations. As can be observed from the construction of $A_\mathcal{S}$ in \eqref{Kaczmarz.Tanabe.iteration.orig}, there are still many inconveniences to use because each column is the product of $Q_i$ and $a_i$.

In this section, we will analyze the inherent structure of the Kaczmarz-Tanabe's iteration \eqref{Kaczmarz.Tanabe.iteration.orig} and derive a concise iterative form similar to the SIRT methods. First, we give the following definitions.

\begin{defn}\label{definition.S_sp} We call $Q_j$ a \textbf{sequential projection matrix} on $a_{j+1},\ldots,a_m$, and denote the sequential projection matrix set with $S_{sp}(a_1,\ldots,a_m)$, i.e.,
\begin{align}\label{sequential.projection.matrix.set}
  S_{sp}(a_1,a_2,\ldots,a_m)=\{Q_1,Q_2,\ldots,Q_{m-1}\}.
\end{align}
\end{defn}

\begin{defn}\label{definition.compatible} For any $Q_i\in S_{sp}$, if there exist $\zeta_{i,1},\ldots, \zeta_{i,m}$ such that
\begin{align}\label{sequential.compatible}
  Q_ia_i=\zeta_{i,1}a_1+\ldots \zeta_{i,m}a_m,
\end{align}
then we call $A$ and $S_{sp}$ \textbf{sequentially compatible}. In general, for any $1\le i\le m, i\le j\le m$, if there exist $\zeta_1^{(i,j)},\ldots,\zeta_m^{(i,j)}$ such that
\begin{align}\label{forward.sequential.compatible}
  Q_ja_i=\zeta_1^{(i,j)}a_1+\ldots \zeta_m^{(i,j)}a_m,
\end{align}
then we call $A$ and $S_{sp}$ \textbf{forward sequentially compatible}, and call $(\zeta_1^{(i,j)},\ldots,\zeta_m^{(i,j)})$ \textbf{compatible vector} of $Q_ja_i$ on $A$.
\end{defn}

\begin{rem}
From Definition \ref{definition.compatible}, we know that, if $A$ and $S_{sp}$ are sequentially compatible, then
\begin{align*}
  a_i^TQ_i^T=\zeta_{i,1}a_1^T+\ldots \zeta_{i,m}a_m^T
\end{align*}
holds, and if $A$ and $S_{sp}$ are forward sequential compatible, then
\begin{align*}
  a_i^TQ_j^T=\zeta_1^{(i,j)}a_1^T+\ldots \zeta_m^{(i,j)}a_m^T
\end{align*}
holds. Therefore, the definitions given in \eqref{sequential.compatible} and \eqref{forward.sequential.compatible} are equivalent to the definition given by their transposes.
\end{rem}

\begin{rem} The definition of forward sequential compatible is actually the constraints on
   \begin{align*}
      a_1^TQ_1^T,\ldots,a_1^TQ_m^T,\quad a_2^TQ_2^T,\ldots,a_2^TQ_m^T,\quad \ldots,\quad a_m^TQ_m^T.
   \end{align*}
Moreover, this definition can be extended completely, but we will not do this because it is beyond the requirements of this paper.
\end{rem}

\begin{rem}
Obviously, $Q_ma_m=a_m$, i.e., $\zeta_m^{(m,m)}\equiv 1$, so the definition of forward sequential compatible can be extended to the case of $i=m$.
\end{rem}

The following theorem shows that $A$ and $S_{sp}$ defined by Kaczmarz's iteration is forward sequential compatible.
\begin{thm}\label{thm.forward.sequential.compatible} Suppose $A$ has no zero row, and $S_{sp}$ is defined by \eqref{sequential.projection.matrix.set}, then $A$ and $S_{sp}$ are forward sequential compatible.
\end{thm}

\proof We take the subscript $(i,j)$ of $Q_ja_i$ as an ordered array and prove the conclusion by mathematical induction.

It is obvious that $a_m^TQ_m^T=a_m^T$, that is, \eqref{forward.sequential.compatible} holds for $(i,j)=(m,m)$ and $(\zeta_1^{(m,m)},\ldots,\zeta_m^{(m,m)})=(0,\ldots,0,1)$. In fact, for any $1\le i\le m$, \eqref{forward.sequential.compatible} obviously holds for $a_i^TQ_m^T$ because $Q_m=I$. Consequently, as the first step of induction, we prove that \eqref{forward.sequential.compatible} holds for $(i,j)=(m-1,m-1)$.
Actually,
\begin{align*}
  a_{m-1}^TQ_{m-1}^T=a_{m-1}^TP_m^TQ_m^T=a_{m-1}^T-\frac{a_{m-1}^T a_m}{\|a_m\|_2^2}a_m^T.
\end{align*}
Hence, \eqref{sequential.compatible} holds for $i=m-1$, where $(\zeta_1^{(m-1,m-1)},\ldots,\zeta_m^{(m-1,m-1)})=(0,\ldots,0,1,-a_{m-1}^Ta_m/\|a_m\|_2^2)$.

Secondly, we suppose \eqref{forward.sequential.compatible} holds for any $(i,j)$ satisfying $s<i<m$ and $s\le t<j<m$, i.e.,
there exists $(\zeta_1^{(i,j)},\ldots,\zeta_m^{(i,j)})$ such that
\begin{align*}
  a_i^TQ_j^T=\zeta_1^{(i,j)}a_1^T+\ldots \zeta_m^{(i,j)}a_m^T.
\end{align*}

Thirdly, we prove that \eqref{forward.sequential.compatible} holds for $(i,j)=(s,t)$. Because of $Q_t=Q_{t+1}P_{t+1}$, then
\begin{align}\label{math.induction.third}
  a_s^TQ_t^T=a_s^TP_{t+1}^TQ_{t+1}^T=a_s^TQ_{t+1}^T-\frac{a_s^Ta_{t+1}}{\|a_{t+1}\|_2^2}a_{t+1}^TQ_{t+1}^T.
\end{align}
From the hypothesis, there exist $(\zeta_1^{(s,t+1)},\ldots,\zeta_m^{(s,t+1)})$ and $(\zeta_1^{(t+1,t+1)},\ldots,\zeta_m^{(t+1,t+1)})$ such that
\begin{align*}
  &a_s^TQ_{t+1}^T=\zeta_1^{(s,t+1)}a_1^T+\ldots+\zeta_m^{(s,t+1)}a_m^T,\\
  &a_{t+1}^TQ_{t+1}^T=\zeta_1^{(t+1,t+1)}a_1^T+\ldots+\zeta_m^{(t+1,t+1)}a_m^T.
\end{align*}
Then, it follows from \eqref{math.induction.third} that
\begin{align*}
  a_s^TQ_t^T&=\zeta_1^{(s,t+1)}a_1^T+\ldots+\zeta_m^{(s,t+1)}a_m^T-\frac{a_s^Ta_{t+1}}{\|a_{t+1}\|_2^2}(\zeta_1^{(t+1,t+1)}a_1^T+\ldots+\zeta_m^{(t+1,t+1)}a_m^T)\\
            &=(\zeta_1^{(s,t+1)}-\frac{a_s^Ta_{t+1}}{\|a_{t+1}\|_2^2}\zeta_1^{(t+1,t+1)})a_1^T+\ldots+(\zeta_m^{(s,t+1)}-\frac{a_s^Ta_{t+1}}{\|a_{t+1}\|_2^2}\zeta_m^{(t+1,t+1)})a_m^T.
\end{align*}
Denote
\begin{align*}
  (\zeta_1^{(s,t)},\ldots,\zeta_m^{(s,t)})=(\zeta_1^{(s,t+1)}-\frac{a_s^Ta_{t+1}}{\|a_{t+1}\|_2^2}\zeta_1^{(t+1,t+1)},\ldots,\zeta_m^{(s,t+1)}-\frac{a_s^Ta_{t+1}}{\|a_{t+1}\|_2^2}\zeta_m^{(t+1,t+1)}).
\end{align*}
This proves that \eqref{forward.sequential.compatible} holds for $(i,j)=(s,t)$.

To sum up the above, the conclusion is proved for all $(i,j)$ with respect to $1\le i\le m,i\le j\le m$. Namely, $A$ and $S_{sp}$ generated by the Kaczmarz's iteration are forward sequentially compatible.  \qed

From Theorem \ref{thm.forward.sequential.compatible}, we have the following decomposition corollary of $A_\mathcal{S}$.
\begin{cor}\label{cor.compatible.decomposition} Under the condition of Theorem \ref{thm.forward.sequential.compatible}, there exists a unit upper triangular matrix $C\in \mathbb{R}^{m\times m}$ such that
\begin{align}
  A_\mathcal{S}=CA.
\end{align}
Here, we call $C$ the \textbf{compatible matrix} of $A$ and $S_{sp}$.
\end{cor}
\proof According to $A_\mathcal{S}=(Q_1 a_1,\ldots,Q_m a_m)^T$ and Theorem \ref{thm.forward.sequential.compatible}, the corollary can be proved by taking $C(i,j)=\zeta_{i,j}$. \qed
\begin{rem}
Corollary \ref{cor.compatible.decomposition} is valuable for the analysis of the Kaczmarz-Tanabe method, which can lead to the standard form of Kaczmarz-Tanabe's iteration (i.e., the matrix-vector form). In fact, it follows from \eqref{Kaczmarz.Tanabe.iteration.orig} and Corollary \ref{cor.compatible.decomposition} that
\begin{align}\label{Kaczmarz.Tanabe.iteration.standard}
  y_{k+1}=y_k+A^TC^TM(b-Ay_k), \qquad k=0,1,2,\ldots.
\end{align}
We can hardly see the shadow of the Kaczmarz iteration from \eqref{Kaczmarz.Tanabe.iteration.standard}, and it is more like a member of SIRT methods. The Kaczmarz's method is known as the algebraic reconstruction technique (ART), However, the appearance of \eqref{Kaczmarz.Tanabe.iteration.standard} makes the boundaries between the ART and SIRT methods confusing, and makes the Kaczmarz method as easy to use as SIRT methods after obtaining $C$.
\end{rem}

In the above, the matrix $C$ exists in theory. For the purpose of dealing with its computational problem, the intuitive idea is to find a matrix $C\in \mathbb{R}^{m\times m}$ that satisfies $A_\mathcal{S}=CA$. For simplicity, we introduce the following notation,
\begin{align}\label{matrix.H}
  H=(h_{i,j}):=AA^TM,
\end{align}
which yields $h_{i,j}=a_i^Ta_j/\|a_j\|_2^2$. For the convenience of description, we introduce the concept of index set.
\begin{defn}\label{definition.index.set} The \textbf{index set} $I_d(n_1,n_2,v)$ is defined as follows
\begin{align*}
  I_d(n_1,n_2,v)=\Big\{[I_d(1),\ldots,I_d(v)]\Big\},
\end{align*}
where $n_1,n_2,v$ are positive integers satisfying $|n_1-n_2|\ge v\ge 2$. $I_d(i)$ is an integer between $n_1$ and $n_2$, and
$I_d(1)=n_1,I_d(v)=n_2$. For any $i<j$, the following is satisfied
\begin{align*}
\begin{array}{ll}
  I_d(i)<I_d(j), & n_1<n_2,\\
  I_d(i)>I_d(j), & n_1>n_2.
\end{array}
\end{align*}
\end{defn}

By the above definition, we know that $I_d(n_1,n_2,v)$ is actually a set of arrays and the elements in every array are arranged by order, e.g.,
\begin{align*}
   &I_d(1,4,2)=\{[1,4]\}, \quad I_d(4,1,2)=\{[4,1]\},\\
   &I_d(1,4,3)=\{[1,2,4],[1,3,4]\},\quad I_d(4,1,3)=\{[4,2,1],[4,3,1]\}.
\end{align*}
We must pay attention to the difference of order. In $[1,4]$, $I_d(1)=1,I_d(2)=4$; and in $[4,1]$, $I_d(1)=4,I_d(2)=1$.

Based on the above definition, we give the expression of $a_i^TQ_i^T\tilde{x}$ when $\tilde{x}\in N(A)^\bot$.
\begin{lem}\label{lemma.compatible.vector} Suppose $A$ has no zero row, $Q_i$ is the sequential projection matrix of $A$ and $\tilde{x}\in N(A)^\bot$. For any $1\le i\le m, i+1\le j\le m-1$, denote
\begin{align}\label{C.element}
   d_{i,j}=\sum_{v=2}^{j-i+1}(-1)^{v-1}\sum_{I_d(i,j,v)}\prod_{s=1}^{v-1} h_{I_d(s),I_d(s+1)}.
\end{align}
Then,
\begin{align}\label{compatible.expression.1}
  a_i^TQ_i^T\tilde{x}=(1,d_{i,i+1},\ldots,d_{i,m})(a_i^T,a_{i+1}^T,\ldots,a_m^T)^T\tilde{x}
\end{align}
holds. That is, the compatible vector of $a_i^TQ_i^T\tilde{x}$ on $A\tilde{x}$ is $(0,\ldots, 0, 1,d_{i,i+1},\ldots,d_{i,m})$.
\end{lem}
\proof When $1\le i\le m$ and $\tilde{x}\in N(A)^\bot$, we have
\begin{align*}
  a_i^TQ_i^T\tilde{x}&=(1,-h_{i,i+1})(a_i^TQ_{i+1}^T\tilde{x},a_{i+1}^TQ_{i+1}^T\tilde{x})^T\nonumber\\
                     &=(1,-h_{i,i+1},-h_{i,i+2}+h_{i,i+1}h_{i+1,i+2})(a_i^TQ_{i+2}^T\tilde{x},a_{i+1}^TQ_{i+2}^T\tilde{x},a_{i+2}^TQ_{i+2}^T\tilde{x})^T\nonumber\\
                     &=(1,-h_{i,i+1},\ldots,\sum_{v=2}^{m-i+1}(-1)^{v-1}\sum_{I_d(i,m,v)}\prod_{s=1}^{v-1} h_{I_d(s),I_d(s+1)})(a_i^TQ_m^T\tilde{x},a_{i+1}^TQ_m^T\tilde{x},\ldots,a_m^TQ_m^T\tilde{x})^T.
\end{align*}
Thus \eqref{compatible.expression.1} holds by taking $d_{i,j}$ according to \eqref{C.element}.\qed

From the proof of Lemma \ref{lemma.compatible.vector}, $d_{i,j}$ is equivalent to the lengthy but intuitive form, i.e.,
\begin{align*}
  \sum_{v=2}^{j-i+1}(-1)^{v-1}\sum_{I_d(i,j,v)}\prod_{s=1}^{v-1} h_{I_d(s),I_d(s+1)}=-h_{i,i+2}+h_{i,i+1}h_{i+1,i+2}+\ldots+(-1)^{j-i}h_{i,i+1}h_{i+1,i+2}\ldots h_{j-1,j}.
\end{align*}

\begin{thm}\label{thm.decomposition.detailed} Under Lemma \ref{lemma.compatible.vector}, let $\Omega=(\omega_{i,j})_{m\times m}$ satisfy
\begin{align}\label{definition.Omega}
  \omega_{i,j}=\left \{
      \begin{array}{ll}
        d_{i,j}, & j>i,\\
        1,       & j=i,\\
        0      , & j<i.
      \end{array}
      \right.
\end{align}
Then,
\begin{align*}
  A_\mathcal{S}=\Omega A
\end{align*}
holds, where $A_\mathcal{S}$ is defined by \eqref{definition.AS}.
\end{thm}
\proof For any $\tilde{x}\in N(A)^\bot$, it follows from \eqref{definition.AS} that
\begin{align}\label{As.x}
  A_\mathcal{S}\tilde{x}=(a_1^TQ_1^T\tilde{x},a_2^TQ_2^T\tilde{x},\ldots,a_m^TQ_m^T\tilde{x})^T.
\end{align}
From Lemma \ref{lemma.compatible.vector} and \eqref{definition.Omega}, we obtain
\begin{align}\label{A.s.equation}
  A_\mathcal{S}\tilde{x}=\Omega A\tilde{x}.
\end{align}
When $\tilde{x}\in N(A)$, \eqref{A.s.equation} obviously holds. Therefore, for any $\tilde{x}\in \mathbb{R}^n$, $A_\mathcal{S}\tilde{x}=\Omega A\tilde{x}$ holds, which means $A_\mathcal{S}=\Omega A$. \qed

Lemma \ref{lemma.compatible.vector} and Theorem \ref{thm.decomposition.detailed} actually show us a specific form of matrix C, i.e., $C\equiv \Omega$, thus we get
\begin{align}\label{Kaczmarz.Tanabe.iteration.specific}
  y_{k+1}=y_k+A^T\Omega^TM(b-Ay_k), \qquad k=0,1,2,\ldots.
\end{align}
If $A$ is a full row rank matrix, the decomposition of $A_{\mathcal{S}}$ is unique.

We specifically refer to \eqref{Kaczmarz.Tanabe.iteration.specific} as the\textbf{ standard form of Kaczmarz-Tanabe's iteration} and still denote by \eqref{Kaczmarz.Tanabe.iteration.standard} with $C=\Omega$.

Let $E(j,i(-h_{j,i}))$ be a matrix obtained by multiplying the $i$-th row of the identity matrix by $-h_{j,i}$ and adding it to the $j$-th row, i.e., the diagonal elements of $E(j,i(-h_{j,i}))$ are all $1$, the $(j,i)$- element is $-h_{j,i}$, and all other elements are $0$.
Consequently, we have the following theorem.

\begin{thm}\label{thm.C.expression} If $\Omega$ is defined as \eqref{definition.Omega}, then
\begin{align}\label{Omega.decomposition.1}
  \Omega=H_1H_2\cdots H_m
\end{align}
holds, where $H_1=I$ and $H_i=\prod\limits_{j=1}^{i-1}E(j,i(-h_{j,i}))$ for any $1<i\le m$.
\end{thm}
\proof For any $\tilde{x}\in N(A)^{\bot}$, we denote $\tilde{b}=A\tilde{x}$. From \eqref{As.x}, we have
\begin{flalign}\label{formula.media.1}
  a_j^TQ_j^T\tilde{x}
  &=(0,\ldots,1,-h_{j,j+1}, \ldots,-h_{j,m}+h_{j,m-1}h_{m-1,m}+\ldots+(-1)^{m-1}h_{j,j+1}h_{j+1,j+2}\cdot\ldots\cdot h_{m-1,m})\nonumber\\
  &\quad\cdot(\tilde{b}_1,\ldots,\tilde{b}_{m-2},\tilde{b}_{m-1},\tilde{b}_m)^T.&
\end{flalign}
From \eqref{As.x}, the coefficient of $\tilde{b}_i (i>j)$ in \eqref{formula.media.1} is actually the $(j,i)$-element of $\Omega$, i.e.,
\begin{align*}
   \omega_{j,i}=-h_{j,i}+h_{j,i-1}h_{i-1,i}+\ldots+(-1)^{i-j}h_{j,j+1}h_{j+1,j+2}\cdot\ldots\cdot h_{i-1,i}.
\end{align*}
Denote $\widehat{H}=H_1\cdots H_m$. In order to show $\Omega=H_1\cdots H_m$, we only need to prove $\omega_{j,i}=\widehat{H}_{j,i}$ for any $i>j$, i.e.,
\begin{align*}
  \omega_{j,i}=e_j^T\widehat{H}e_i,
\end{align*}
where $e_j$ and $e_i$ are the $j$th and $i$th columns of the identity matrix in $\mathbb{R}^{m\times m}$ \cite[p72]{David_linear_algebra}, respectively. Owing to $e_j^TH_k=e_j^T$ when $j\ge k$, and $H_le_i=e_i$ when $i\neq l$, it follows that when $i>j$,
\begin{align*}
  e_j^T\widehat{H}e_i&=e_j^TH_{j+1}\cdots H_ie_i\nonumber\\
                     &=(e_j^TH_{j+1})H_{j+2}\cdots H_ie_i\nonumber\\
                     &=((0,\ldots,1,-h_{j,j+1},0,\ldots,0)H_{j+2})H_{j+3}\cdots H_ie_i\nonumber\\
       &=(0,\ldots,1,-h_{j,j+1},\ldots,-h_{j,i}+h_{j,i-1}h_{i-1,i}+\ldots+(-1)^{i-j}h_{j,j+1}\cdots h_{i-1,i},0,\ldots,0)e_i\nonumber\\
       &=-h_{j,i}+h_{j,i-1}h_{i-1,i}+\ldots+(-1)^{i-j}h_{j,j+1}h_{j+1,j+2}\cdots h_{i-1,i}.&
\end{align*}
This proves $\omega_{j,i}=\widehat{H}_{j,i}$ for any $1\le j\le n-1$ and $i>j$. Additionally, $\omega_{j,j}=\widehat{H}_{j,j}=1$ holds for any $1\le j\le n$. Consequently, the conclusion is proved. \qed

Theorem \ref{thm.C.expression} actually gives the calculation formula of $\Omega$ defined in \eqref{definition.Omega}. However, it is not a good idea to calculate $\Omega$ directly according to \eqref{Omega.decomposition.1} because the calculation speed may be slow. In fact, the matrix $\Omega$ can be calculated in parallel mode by dividing the multiplication of $H_1H_2\cdots H_m$ into several small parts, but we should notice that the block operation is executed on matrix $\Omega$ but not on the whole linear system. Consequently, performing the block operation on linear system and solving each linear subsystem with the Kaczmarz-Tanabe method, which indeed can reduce the cost of calculating $\Omega$, will derive the block Kaczmarz-Tanabe method.

\section{The standard form of symmetric Kaczmarz-Tanabe method and its convergence}
\label{section.the.symmetric.Kacmarz.Tanbe}

In this section, we mainly consider the standard form of the symmetric Kaczmarz-Tanabe's iteration and analyze its convergence rate, and then compare it with the convergence rate of the Kaczmarz-Tanabe's iteration.

Let $\{x_k, k>0\}$ be the vector sequence determined by \eqref{symmetric.Kaczmarz.orig} and \eqref{symmetric.Kaczmarz.con}. Denote
\begin{align}\label{period.iteration}
  \bar{y}_{k+1}=x_{k\cdot(2m-2)+m},\quad y_{k+1}=x_{(k+1)\cdot(2m-2)}, \quad k=0,1,\ldots.
\end{align}
Then, from \eqref{Kaczmarz.Tanabe.iteration.standard},
\begin{align}\label{the.symmetric.Kaczmarz.iteration}
  \bar{y}_{k+1}=y_k+A^TC^TM(b-Ay_k)
\end{align}
holds, which is indeed the Kaczmarz-Tanabe's iteration from $y_k$ to $\bar{y}_{k+1}$.

Next, we consider the iterative formula of Kaczmarz-Tanabe method for Kaczmarz's projection from equation $m-1$ to equation $2$ in reverse order, i.e., the Kaczmarz-Tanabe's iteration from $\bar{y}_{k+1}$ to $y_{k+1}$. Define
\begin{align}\label{definition.barQ.i}
   \bar{Q}_i=P_2\ldots P_{i-1},\quad i=3,\ldots,m-1.
\end{align}
Thus, $\bar{Q_i}$ is the \textbf{sequential projection matrix} on $(a_{i-1},\ldots,a_2)^T$, and
\begin{align}\label{property.bar.Qi}
  \bar{Q}_i=\bar{Q}_{i-1}P_{i-1}, \quad i=3,\ldots,m-1.
\end{align}
The sequential projection matrix set reads
\begin{align*}
  \bar{S}_{sp}(a_{m-1},\ldots,a_2)=\{\bar{Q}_3,\ldots,\bar{Q}_{m-1}\}.
\end{align*}
Additionally, we have
\begin{align}\label{appointment.2}
  \bar{Q}_1=\bar{Q}_m={\bf 0}\in\mathbb{R}^{m\times m}, \quad \bar{Q}_2=I,
\end{align}
and denote
\begin{align}\label{agreement.barQ}
  \bar{Q}:=P_2\ldots P_{m-1}.
\end{align}
Hence the Kaczmarz's projections from equation $m-1$ to $2$ are equivalent to
\begin{align}\label{the.quasi.symmetric.Kaczmarz}
  y_{k+1}=\bar{Q}\bar{y}_{k+1}+\bar{A}_s^TMb,
\end{align}
where
\begin{align}\label{definition.bar.AS}
  \bar{A}_\mathcal{S}=(\bar{Q}_1a_1,\bar{Q}_2a_2,\bar{Q}_3a_3,\ldots,\bar{Q}_{m-1}a_{m-1},\bar{Q}_ma_m)^T.
\end{align}

Note that \eqref{the.quasi.symmetric.Kaczmarz} is not the symmetric Kaczmarz-Tanabe iteration but the symmetric part of the symmetric Kaczmarz's iteration, i.e., the case of Kaczmarz's projection \eqref{symmetric.Kaczmarz.orig} for $i=m-1,\ldots, 2$.

Before deriving the standard form of the symmetric Kaczmarz-Tanabe's iteration, we first give the relationship between $\bar{Q}$ and $\bar{A}_\mathcal{S}$ appearing in \eqref{the.quasi.symmetric.Kaczmarz}.

\begin{lem}\label{lemma.relationship.barQ.barAs} Suppose $A$ has no zero row, and $\bar{Q}$ and $\bar{A}_\mathcal{S}$ are defined as \eqref{agreement.barQ} and \eqref{definition.bar.AS}. Then,
\begin{align}\label{relationship.barQ.AsTMa}
  \bar{Q}=I-\bar{A}_\mathcal{S}^TMA,
\end{align}
where $M$ is defined in \eqref{definition.M}.
\end{lem}
\proof From \eqref{agreement.barQ} and \eqref{property.bar.Qi}, we have
\begin{align*}
  \bar{Q}&=P_2\ldots P_{m-1}\\
         &=\bar{Q}_{m-1}-\bar{Q}_{m-1}\frac{a_{m-1}a_{m-1}^T}{\|a_{m-1}\|_2^2}\\
         &= \cdots\\
         &=\bar{Q}_2-\bar{Q}_2\frac{a_2a_2^T}{\|a_2\|_2^2}-\bar{Q}_3\frac{a_3a_3^T}{\|a_3\|_2^2}-\ldots-\bar{Q}_{m-2}\frac{a_{m-2}a_{m-2}^T}{\|a_{m-2}\|_2^2}-\bar{Q}_{m-1}\frac{a_{m-1}a_{m-1}^T}{\|a_{m-1}\|_2^2}-\bar{Q}_m\frac{a_ma_m^T}{\|a_m\|_2^2}.
\end{align*}
From \eqref{appointment.2}, it follows
\begin{align*}
         \bar{Q}&=I-\bar{Q}_1\frac{a_1a_1^T}{\|a_1\|_2^2}-\bar{Q}_2\frac{a_2a_2^T}{\|a_2\|_2^2}-\bar{Q}_3\frac{a_3a_3^T}{\|a_3\|_2^2}-\ldots-\bar{Q}_{m-2}\frac{a_{m-2}a_{m-2}^T}{\|a_{m-2}\|_2^2}-\bar{Q}_{m-1}\frac{a_{m-1}a_{m-1}^T}{\|a_{m-1}\|_2^2}-\bar{Q}_m\frac{a_ma_m^T}{\|a_m\|_2^2}\\
         &=I-(\bar{Q}_1a_1,\bar{Q}_2a_2,\ldots,\bar{Q}_ma_m)\text{diag}(\frac{1}{\|a_1\|_2^2},\ldots,\frac{1}{\|a_m\|_2^2})(a_1,a_2,\ldots,a_m)^T.
\end{align*}
This proves \eqref{relationship.barQ.AsTMa}.\qed

According to Lemma \ref{lemma.relationship.barQ.barAs}, we get the equivalent form of \eqref{the.quasi.symmetric.Kaczmarz},
\begin{align}\label{the.symmetric.Kaczmarz.4}
  y_{k+1}=\bar{y}_{k+1}+\bar{A}_\mathcal{S}^TM(b-A\bar{y}_{k+1}).
\end{align}
We should notice that \eqref{the.symmetric.Kaczmarz.4} is not the final form of the symmetric Kaczmarz-Tanabe's iteration because it does not include the Kaczmarz projection process from $i=1$ to $m$. We next consider the matrix-vector form of \eqref{the.symmetric.Kaczmarz.4}. First, we have the following existence theorem.
\begin{thm}\label{thm.decomposition.bar.As} Suppose $A$ has no zero row, then there exists $\widehat{C}$ such that
\begin{align}
    \bar{A}_\mathcal{S}=\widehat{C}A.
\end{align}
\end{thm}
\proof Similar to Theorem \ref{thm.forward.sequential.compatible} and Corollary \ref{cor.compatible.decomposition}, the existence of $\widehat{C}$ can be proved. We omit the process here.\qed

Because $\bar{Q}_1=\bar{Q}_m={\bf 0}$, the elements in the first and the last rows of $\widehat{C}$ are zero. According to \eqref{definition.bar.AS}, $\bar{A}_\mathcal{S}$ has nothing to do with $a_1$ when $\bar{Q}_1={\bf 0}$, which implies that the first column of $\widehat{C}$ is zero vector. These characteristics are the major difference between $\widehat{C}$ and $C$. Before considering the specific expression of $\widehat{C}$, we first introduce the following lemma.

\begin{lem}\label{lemma.compatible.barQ} Suppose $A$ has no zero row and $\bar{Q}_i$ is defined as \eqref{definition.barQ.i}. For $3\le i\le m-1, 2\le j\le i-1$, denote
\begin{align}\label{barC.element}
   \bar{d}_{i,j}=\sum_{v=2}^{i-j+1}(-1)^{v-1}\sum_{I_d(i,j,v)}\prod_{s=1}^{v-1} h_{I_d(s),I_d(s+1)},
\end{align}
where $h_{\centerdot,\centerdot}$ is defined by \eqref{matrix.H}. Then, for any $\bar{x}\in N(A)^\bot$,
\begin{align}\label{compatible.expression}
  a_i^T\bar{Q}_i^T\bar{x}=(0,\bar{d}_{i,2},\ldots,\bar{d}_{i,i-1},1,\ldots,0)(a_1,\ldots,a_m)^T\bar{x}
\end{align}
holds.
\end{lem}
\proof Obviously, when $3\le i\le m-1$ and $\bar{x}\in N(A)^\bot$, from \eqref{property.bar.Qi} we have
\begin{align}\label{compatible.barQ}
  a_i^T\bar{Q}_i^T\bar{x}&=(-h_{i,i-1},1)(a_{i-1}^T\bar{Q}_{i-1}^T\bar{x},a_i^T\bar{Q}_{i-1}^T\bar{x})^T\nonumber\\
  &=(-h_{i,i-1},1)(a_{i-1}^T\bar{Q}_{i-1}^T\bar{x},a_i^T\bar{Q}_{i-1}^T\bar{x})^T\nonumber\\
  &=(-h_{i,i-2}+h_{i,i-1}h_{i-1,i-2},-h_{i,i-1},1)(a_{i-2}^TQ_{i-2}^T\bar{x},a_{i-1}^TQ_{i-2}^T\bar{x},a_i^TQ_{i-2}^T\bar{x})^T\nonumber\\
  &=\ldots\nonumber\\
  &=(-h_{i,2}+\sum_{k=3}^{i-1}h_{i,k}h_{k,2}+\ldots +(-1)^{i-2}\prod_{k=2}^{i-1}h_{k+1,k},\ldots,-h_{i,i-1},1)(a_2^T\bar{x},\ldots,a_i^T\bar{x})^T.
\end{align}
Taking $\bar{d}_{i,j}$ in \eqref{compatible.barQ} according to \eqref{barC.element} yields
\begin{align*}
  a_i^T\bar{Q}_i^T\bar{x}=(\bar{d}_{i,2},\ldots,\bar{d}_{i,i-1},1)(a_2,\ldots,a_i)^T\bar{x}.
\end{align*}
This proves \eqref{compatible.expression}.\qed

Similar to Theorem \ref{thm.decomposition.detailed}, we have the following theorem.
\begin{thm}\label{thm.decomposition.bar.As.2} Under the condition of Lemma \ref{lemma.compatible.barQ}, let $\widehat\Omega=(\widehat\omega_{i,j})_{m\times m}$ satisfy
\begin{align}\label{definition.bar.Omega}
  \widehat\omega_{i,j}=\left \{
      \begin{array}{ll}
        \bar{d}_{i,j}, & 1<i<m,2<j<i-1,\\
        0,             & 1<i<m,j=1 \vee j>i,\\
        1,             & 1<i<m,j=i,\\
        0,             & i=1 \vee m, 1\le j\le m.
      \end{array}
      \right.
\end{align}
Then,
\begin{align}
  \bar{A}_\mathcal{S}=\widehat\Omega A
\end{align}
holds, where $\bar{A}_\mathcal{S}$ is defined as \eqref{definition.bar.AS}.
\end{thm}
\proof For any $\bar{x}\in N(A)^\bot$, it follows from \eqref{definition.bar.AS} that
\begin{align}\label{bar.As.x}
  \bar{A}_\mathcal{S}\bar{x}=(a_1^T\bar{Q}_1^T\bar{x},a_2^T\bar{Q}_2^T\bar{x},\ldots,a_m^T\bar{Q}_m^T\bar{x})^T.
\end{align}
By Lemma \ref{lemma.compatible.barQ} and \eqref{definition.bar.Omega}, then we get
\begin{align}\label{barA.s.equation}
  \bar{A}_\mathcal{S}\bar{x}=\widehat\Omega A\bar{x}.
\end{align}
When $\bar{x}\in N(A)$, from \cite[Corollary 2.2]{Kang2021Convergence}, $\bar{A}_\mathcal{S}\bar{x}=0$ holds. Thus, \eqref{barA.s.equation} also holds. Then, for any $\bar{x}\in \mathbb{R}^n$, $\bar{A}_\mathcal{S}\bar{x}=\widehat\Omega A\bar{x}$ holds, which means $\bar{A}_\mathcal{S}=\widehat\Omega A$. \qed

Let
\begin{align}
  \hat{E}(j,i(-h_{j,i}))(s,t)=\left \{
  \begin{array}{ll}
     E(j,i(-h_{j,i}))(s,t), &\quad (s,t)\neq (1,1)\wedge (m,m),\\
     0,                       &\quad (s,t)=(0,0)\vee (m,m).
  \end{array}
  \right.
\end{align}
Similar to Theorem \ref{thm.C.expression}, we have the following decomposition of $\hat{\Omega}$.
\begin{thm}\label{thm.hat.Omega.expression} If $\widehat\Omega$ is defined as in Theorem \ref{thm.decomposition.bar.As.2}, then
\begin{align}\label{Omega.decomposition}
  \widehat\Omega=\hat{H}_{m-1}\hat{H}_{m-2}\cdots \hat{H}_2
\end{align}
holds, where $\hat{H}_i=\prod\limits_{j=i+1}^{m-1}\hat{E}(j,i(-h_{j,i}))$ for any $2\le i\le m-1$.
\end{thm}
\proof Denote $\widetilde{H}=\hat{H}_{m-1}\hat{H}_{m-2}\cdots \hat{H}_2$. Obviously, $\widehat\Omega$ and $\widetilde{H}$ are unit lower triangular matrices with the same order, so we only need to prove that the non-zero elements are equal. For any $\bar{x}\in N(A)^{\bot}$, from \eqref{compatible.barQ} and $\bar{Q}_2=I$, we have
\begin{align}\label{formula.media.2}
  a_i^T\bar{Q}_i^T\bar{x}&=(-h_{i,2}+\sum_{k=3}^{i-1}h_{i,k}h_{k,2}+\ldots +(-1)^{i-2}\prod_{k=2}^{i-1}h_{k+1,k},\ldots,-h_{i,i-1},1)(a_2^T\bar{x},\ldots,a_i^T\bar{x})^T.&
\end{align}
In \eqref{formula.media.2}, the coefficient of $a_j^T\bar{x} (2\le j<i)$ is actually the $(i,j)$-element of $\widehat\Omega$, i.e.,
\begin{align*}
   \hat\omega_{i,j}=-h_{i,j}+\sum\limits_{k=j+1}^{i-1}h_{i,k}h_{k,j}+(-1)^{i-j}\prod_{k=j}^{i-1}h_{k+1,k}.
\end{align*}
In order to show $\widehat\Omega=\hat{H}_{m-1}\hat{H}_{m-2}\cdots \hat{H}_2$, we only need to prove $\hat\omega_{i,j}=\widetilde{H}_{i,j}$ (where $\widetilde{H}_{i,j}$ denotes the $(i,j)$-element of $\widetilde{H}$), i.e.,
\begin{align*}
  \hat\omega_{i,j}=e_i^T\widetilde{H}e_j.
\end{align*}
Owing to $e_i^T\hat{H}_k=e_i^T$ when $i\le k$ and $i=m$, and $\hat{H}_le_j=e_j$ when $j\neq l$, we have for $2\le i\le m-1$ and $2\le j<i$,
\begin{align*}
  e_i^T\widetilde{H}e_j
       &=e_i^T\hat{H}_{i-1}\hat{H}_{i-2}\cdots \hat{H}_je_j\nonumber\\
       &=(e_i^T\hat{H}_{i-1})\hat{H}_{i-2}\cdots \hat{H}_je_j\nonumber\\
       &=((0,\ldots,-h_{i,i-1},1,0,\ldots,0)\hat{H}_{i-2})H_{i-3}\cdots \hat{H}_je_j\nonumber\\
       &=((0,\ldots,0,-h_{i,i-2}+h_{i,i-1}h_{i-1,i-2},-h_{i,i-1},1,0,\ldots,0)\hat{H}_{i-3})\cdots \hat{H}_jej\nonumber\\
       &=-h_{i,j}+\sum_{k=j+1}^{i-1}h_{i,k}h_{k,j}+\ldots+(-1)^{i-j}h_{i,i-1}h_{i-1,i-2}\cdots h_{j+1,j}.&
\end{align*}
This proves $\hat\omega_{i,j}=\widetilde{H}_{i,j}$ for any $2\le i\le m-1$ and $2\le j<i$. Additionally, $\hat{\omega}_{i,i}=\widetilde{H}_{i,i}=1$ holds for any $2\le i\le m-1$. Consequently, the conclusion is proved. \qed

Compared with Theorem \ref{thm.decomposition.bar.As}, $\widehat\Omega$ in Theorem \ref{thm.decomposition.bar.As.2} gives the specific form of $\widehat{C}$ and is still denoted by $\widehat{C}$. Thus, from \eqref{the.symmetric.Kaczmarz.4} we obtain \begin{align}\label{symmetric.Kaczmarz.iteration}
  y_{k+1}=\bar{y}_{k+1}+A^T\widehat{C}^T M(b-A\bar{y}_{k+1}).
\end{align}

Based on \eqref{the.symmetric.Kaczmarz.iteration} and \eqref{symmetric.Kaczmarz.iteration}, we obtain the following theorem.
\begin{thm}\label{thm.symmetric.Kaczmarz.Tanabe} Suppose $A$ has no zero row. Then, there exists matrix $\bar{C}\in \mathbb{R}^{m\times m}$, such that the \textbf{symmetric Kaczmarz-Tanabe's iteration} can be written as
\begin{align}\label{symmetric.Kaczmarz.Tanabe.iteration.specific}
  y_{k+1}=y_k+A^T\bar{C}^TM(b-Ay_k).
\end{align}
\end{thm}
\proof From \eqref{the.symmetric.Kaczmarz.iteration} and \eqref{symmetric.Kaczmarz.iteration}, the symmetric Kaczmarz-Tanabe's iteration is given by
\begin{align*}
  y_{k+1}&=\bar{y}_{k+1}+A^T\widehat{C}^TM(b-A\bar{y}_{k+1})\\
         &=y_k+A^T(\widehat{C}^T+C^T-\widehat{C}^TMAA^TC^T)M(b-Ay_k).
\end{align*}
Denote $\bar{C}:=\widehat{C}+C-CAA^TM\widehat{C}$, then \eqref{symmetric.Kaczmarz.Tanabe.iteration.specific} is proved.\qed

From \eqref{Kaczmarz.Tanabe.iteration.specific} and \eqref{symmetric.Kaczmarz.Tanabe.iteration.specific}, we know that the Kaczmarz-Tanabe's iteration and the symmetric Kaczmarz-Tanabe's iteration have the same matrix-vector form. Then, from \eqref{symmetric.Kaczmarz.Tanabe.iteration.specific}, we also have the following equivalent expression
\begin{align*}
  y_{k+1}=(I-A^T\bar{C}^TMA)y_k+A^T\bar{C}^TMb,
\end{align*}
where $I-A^T\bar{C}^TMA$ is the iteration matrix of the symmetric Kaczmarz-Tanabe method. Considering the principle of the symmetric Kaczmarz's iteration, we have the following corollary.
\begin{cor} Suppose $A$ has no zero row. Then, for the symmetric Kaczmarz-Tanabe's iteration,
\begin{align}\label{barQ.Q}
  \bar{Q}Q=P_2\ldots P_{m-1}P_m\ldots P_1=I-A^T\bar{C}^TMA
\end{align}
holds, where $\bar{C}$ is consistent with that in Theorem \ref{thm.symmetric.Kaczmarz.Tanabe}.
\end{cor}

Let $e_k=y_k-P_{N(A)}y_0-x^\dagger$. Then, it follows from \eqref{symmetric.Kaczmarz.Tanabe.iteration.specific} that
\begin{align}\label{error.symm.Kaczmarz.Tanabe}
  e_{k+1}=(I-A^T\bar{C}^TMA)e_k.
\end{align}
For the symmetric Kaczmarz-Tanabe's iteration, the following holds.

\begin{thm}\label{thm.set.property.symm.KT} For any initial vector $y_0\in \mathbb{R}^n$, let $\{y_k,k>0\}$ be generated by the symmetric Kaczmarz-Tanabe's iteration \eqref{symmetric.Kaczmarz.Tanabe.iteration.specific}. Then,
\begin{align*}
  e_k\in N(A)^\bot
\end{align*}
holds.
\end{thm}
\proof We prove the conclusion by mathematical induction. First, the fact $e_0\in N(A)^\bot$ holds because $y_0-P_{N(A)}y_0$ and $x^\dagger$ belong to $N(A)^\bot$. Second, if we assume that for any given $k\ge 0$, $e_k\in N(A)^\bot$, then $e_{k+1}\in N(A)^\bot$, this is because for any $z\in N(A)$,
\begin{align*}
  \langle e_{k+1},z\rangle=\langle (I-A^T\bar{C}^TMA)e_k,z\rangle=\langle e_k,z\rangle-\langle \bar{C}^TMAe_k,Az\rangle=0.
\end{align*}
Which proves the conclusion.\qed

\begin{lem}\label{lemma.orthogonality} For any $1\le i\le m$ and $x\in N(A)^\bot$,
\begin{align*}
  P_ix\in N(A)^\bot
\end{align*}
holds. That is, $N(A)^\bot$ is an invariant subspace for any $P_i$.
\end{lem}
\proof For any $1\le i\le m$ and $z\in N(A)$,
\begin{align*}
  \langle P_ix,z\rangle=\langle (I-\frac{a_ia_i^T}{\|a_i\|_2^2})x,z\rangle=\langle x,z\rangle-\frac{1}{\|a_i\|_2^2}\langle a_i^Tx,a_i^Tz\rangle=0.
\end{align*}
holds.
\qed

By Lemma \ref{lemma.orthogonality}, we can obtain the following estimation of $\|P_1e_{k+1}\|_2$.

\begin{thm}\label{thm.convergence.rate.symm.Kaczmarz.Tanabe.1} Under the condition of Theorem \ref{thm.set.property.symm.KT},
\begin{align}\label{result.convergence.rate.symm.Kaczmarz.Tanabe}
  \|P_1e_{k+1}\|_2\le \max_{0<\sigma_i<1}\sigma_i^2\|P_1e_k\|_2,\quad k=0,1,\ldots
\end{align}
holds, where $\sigma_i$ is a singular value of $Q$.
\end{thm}
\proof From \eqref{error.symm.Kaczmarz.Tanabe}, we have
\begin{align}\label{mapping.error.1}
  P_1e_{k+1}=P_1(I-A^T\bar{C}^TMA)e_k.
\end{align}
Note that $I-A^T\bar{C}^TMA=P_2\ldots P_{m-1}P_m\ldots P_1$, then
\begin{align}\label{matrix.property.1}
   P_1(I-A^T\bar{C}^TMA)=P_1P_2\ldots P_{m-1}P_m\ldots P_1=Q^TQP_1,
\end{align}
thus
\begin{align}\label{mapping.error}
  P_1e_{k+1}=Q^TQP_1e_k.
\end{align}
Moreover, by Lemma \ref{lemma.orthogonality}, $P_1e_k\in N(A)^\bot$. From Theorem \ref{thm.set.property.symm.KT} and \cite[Theorem 1.3]{Kang2021Convergence}, we have $\|Q|_{N(A)^\bot}\|_2<1$, then we get
\begin{align*}
  \|P_1e_{k+1}\|_2\le \max_{0<\sigma_i<1}\sigma_i^2\|P_1e_k\|_2,
\end{align*}
where $\sigma_i$ is a singular value of $Q$ .
\qed
\begin{cor}\label{cor.weak.convergence.rate.symm.Kaczmarz.Tanabe} Under the condition of Theorem \ref{thm.set.property.symm.KT}, for some $k\ge 0$, if $e_{k+1}\in N(P_1)^\bot$, then
\begin{align}\label{result.convergence.rate.1}
   \|e_{k+1}\|_2\le\max_{0<\sigma_i<1}\sigma_i^2\|P_1^\dagger\|_2\|e_k\|_2
\end{align}
holds, where $\sigma_i$ is the singular value of $Q$ and $P_1^\dagger$ denotes the pseudo-inverse of $P_1$.
\end{cor}
\proof When $e_{k+1}\in N(P_1)^\bot$, we have
\begin{align}\label{key.formula.1}
  P_1^\dagger P_1e_{k+1}=e_{k+1}.
\end{align}
Hence,
\begin{align*}
  \|e_{k+1}\|_2\le \|P_1^\dagger\|_2\|P_1e_{k+1}\|_2,
\end{align*}
and from \eqref{result.convergence.rate.symm.Kaczmarz.Tanabe}, we obtain \eqref{result.convergence.rate.1}.\qed

The equality \eqref{key.formula.1} depends on $e_{k+1}\in N(P_1)^\bot$. If the latter is not satisfied, then \eqref{result.convergence.rate.1} may not hold. In the following theorem, we give a general conclusion without the constraint condition $e_{k+1}\in N(P_1)^\bot$.

\begin{thm}\label{thm.weak.convergence.rate.symm.Kaczmarz.Tanabe} Under the condition of Theorem \ref{thm.set.property.symm.KT}, for any $k\ge 0$, at least one of the following statements is true,

{\normalfont (i)} \quad  $\|e_{k+1}\|_2<\max\limits_{0<\sigma_i<1}\sigma_i\|e_k\|_2$;

{\normalfont(ii)} \quad $\|e_{k+2}\|_2<\max\limits_{0<\sigma_i<1}\sigma_i^2\|e_k\|_2$;

\noindent where $\sigma_i$ is a singular value of $Q$.
\end{thm}

\proof First, we have $N(A)=N(a_1^T)\cap N(a_2^T)\cap\ldots \cap N(a_m^T)$. Then,
\begin{align*}
  N(A)^\bot=N(a_1^T)^\bot\cup N(a_2^T)^\bot\cup\cdots\cup N(a_m^T)^\bot.
\end{align*}
Recall that $e_k\in N(A)^\bot$ and $Qe_k\in N(A)^\bot$, then, at least one of ($I_1$) and ($I_2$) holds:

($I_1$) Among $P_2,\ldots,P_m$, there exists at least one $P_i$ such that
\begin{align}\label{compress.property}
  \|P_iQe_k\|_2<\|Qe_k\|_2.
\end{align}

($I_2$) $\|P_1Qe_k\|_2<\|Qe_k\|_2$.\\
Since $Qe_k\in N(A)^\bot$,  either $Qe_k\in N(a_2^T)^\bot\cup\cdots\cup N(a_m^T)^\bot$ or $Qe_k\in N(a_1)^\bot$. When $Qe_k\in N(a_2^T)^\bot\cup\cdots\cup N(a_m^T)^\bot$, without loss of generality, we suppose $Qe_k\in N(a_m^T)^\bot$. Thus
\begin{align*}
  \|P_mQe_k\|_2^2=\langle Qe_k-\frac{a_ma_m^T}{\|a_m\|_2^2}Qe_k, Qe_k-\frac{a_ma_m^T}{\|a_m\|_2^2}Qe_k\rangle=\|Qe_k\|_2^2-\frac{(a_m^TQe_k)^2}{\|a_m\|_2^2}<\|Qe_k\|_2^2.
\end{align*}
i.e., $\|P_mQe_k\|_2<\|Qe_k\|_2$. If $Qe_k\in N(a_1^T)^\bot$, then $\|P_1Qe_k\|_2<\|Qe_k\|_2$. Consequently, when $Qe_k\in N(A)^\bot$, at least one of ($I_1$) and ($I_2$) holds.

When $(I_1)$ holds, let $l$ be the largest index $i$ that satisfies \eqref{compress.property}, i.e., $P_iQe_k=Qe_k$ for $l+1\le i\le m$. If $Qe_k\in N(A)^\bot$, from Theorem \ref{thm.set.property.symm.KT} and Lemma \ref{lemma.orthogonality}, we have
\begin{align*}
  \|e_{k+1}\|_2&=\|P_2\ldots P_{m-1}P_mQe_k\|_2\le\|P_lQe_k\|_2<\|Qe_k\|_2.
\end{align*}
Therefore, when $e_k\in N(A)^\bot$,
\begin{align*}
  \|Qe_k\|_2\le \max_{0<\sigma_i<1}\sigma_i\|e_k\|_2,
\end{align*}
this proves statement {\normalfont (i)}.

When ($I_2$) holds, we assume that $Qe_k\in N(a_2^T)\cap\ldots \cap N(a_m^T)$, then
\begin{align*}
  e_{k+1}=P_2\cdots P_{m-1}Qe_k=Qe_k
\end{align*}
holds. Moreover,
\begin{align*}
  e_{k+2}=P_2P_3\cdots P_{m-1}P_m\cdots P_2P_1e_{k+1}=P_2P_3\cdots P_{m-1}QP_1Qe_k.
\end{align*}
Then,
\begin{align*}
  \|e_{k+2}\|_2\le \|QP_1Qe_k\|_2
\end{align*}
holds. Since $P_1Qe_k\in N(A)^\bot$,
\begin{align*}
  \|e_{k+2}\|_2\le \max_{0<\sigma_i<1}\sigma_i\|P_1Qe_k\|_2<\max_{0<\sigma_i<1}\sigma_i\|Qe_k\|_2<\max_{0<\sigma_i<1}\sigma_i^2\|e_k\|_2
\end{align*}
hold. This proves statement (ii).\qed

\begin{rem} From Theorem \ref{thm.weak.convergence.rate.symm.Kaczmarz.Tanabe} we can see that the convergence rate of the symmetric Kaczmarz-Tanabe method is better than that of the Kaczmarz-Tanabe method (since `$\le$' is replaced by `$<$'). However, the comparison is actually unfair because each iteration of the symmetric Kaczmarz-Tanabe method performs $2m-2$ orthogonal projections, while the Kaczmarz-Tanabe method only makes $m$ orthogonal projections. Consequently, we'd better compare the convergence rate of the symmetric Kaczmarz-Tanabe method with that of the two-step Kaczmarz-Tanabe method. Supposing $\{y_k,k>0\}$ is the sequence of the Kaczmarz-Tanabe's iteration, so the two-step Kaczmarz-Tanabe's iteration can be represented by $z_k=y_{2k}$. Let $\bar{e}_k=z_k-x^\dagger-P_{N(A)}x_0$, then
\begin{align*}
  \|\bar{e}_{k+1}\|_2\le \max_{0<\sigma_i<1}\sigma_i^2\|\bar{e}_k\|_2.
\end{align*}
According to Theorem \ref{thm.weak.convergence.rate.symm.Kaczmarz.Tanabe}(i), the convergence rate of the two-step Kaczmarz-Tanabe's iteration is better than that of the symmetric Kaczmarz-Tanabe's iteration.
\end{rem}
\begin{rem}
As can be seen from \eqref{Kaczmarz.Tanabe.iteration.standard} and \eqref{symmetric.Kaczmarz.Tanabe.iteration.specific}, the Kaczmarz-Tanabe method and the symmetric Kaczmarz-Tanabe method have the same iterative formula, but $C$ is different from $\bar{C}$.
When $C$ and $\bar{C}$ are known, one iteration of the Kaczmarz-Tanabe method is equivalent to $m$ Kaczmarz's iterations, while one iteration of the symmetric Kaczmarz-Tanabe method is equivalent to $2m-2$ Kaczmarz's iterations. From this point of view, the calculation efficiency of the symmetric Kaczmarz-Tanabe method is higher than that of the Kaczmarz-Tanabe method.
\end{rem}

\section{The related algorithms}
\label{section.algorithm}
For the Kaczmarz-Tanabe method, the core work is to generate matrix $C$. Once $C$ is obtained, the Kaczmarz-Tanabe's iteration is easy to perform. Algorithm \ref{algorithm.C} shows the process flow of calculating $C$.
\begin{breakablealgorithm}\label{algorithm.C}
\caption{The calculation of matrix $C$}\label{algorithm.C.flow}
\begin{algorithmic}[1]
\State \textbf{Input}
\State \indent $A=(a_1,a_2,\ldots,a_m)^T$
\State \indent $C=I_m$ \Comment{$I_m$ is an identity matrix with order $m$}
\State \indent $k \gets m$
\While{$k>1$}
  \State $i \gets k$
  \While{$i>1$}
     \State $j \gets m$
     \While{$j>k-1$}
     \If {$a_k^Ta_k=0$}
        \State $C(i-1,j)=C(i-1,j)$
     \Else
        \State $C(i-1,j)=C(i-1,j)+(-a_{i-1}^Ta_k/a_k^Ta_k)C(k,j)$
     \EndIf
        \State $j \gets j-1$
     \EndWhile
  \State $i \gets i-1$
  \EndWhile
\State $k \gets k-1$
\EndWhile
\State \textbf{Output} $C$
\end{algorithmic}
\end{breakablealgorithm}

For the symmetric Kaczmarz-Tanabe method, $\bar{C}=\widehat{C}^T+C^T-\widehat{C}^TMAA^TC^T$, where $C$ is the matrix obtained by Algorithm \ref{algorithm.C}. Therefore, we only need to compute $\widehat{C}$ in order to perform the symmetric Kaczmarz-Tanabe's iteration. Algorithm \ref{algorithm.C.hat} shows the process flow for computing $\widehat{C}$.

\begin{breakablealgorithm}\label{algorithm.C.hat}
\caption{The calculation of matrix $\widehat{C}$}\label{algorithm.C.hat.flow}
\begin{algorithmic}[1]
\State \textbf{Input}
\State \indent $A=(a_1,a_2,\ldots,a_m)^T$
\State \indent $\widehat{C}=I_m$ \Comment{$I_m$ is an identity matrix with order $m$}
\State \indent $k \gets m-1$
\While{$k>1$}
  \State $i \gets m-1$
  \While{$i>k$}
  j  \State $j \gets k$
     \While{$j>1$}
     \If {$a_k^Ta_k=0$}
        \State $\widehat{C}(i,j)=\widehat{C}(i,j)$
     \Else
        \State $\widehat{C}(i,j)=\widehat{C}(i,j)+(-a_i^Ta_k/a_k^Ta_k)\widehat{C}(k,j)$
     \EndIf
        \State $j \gets j-1$
     \EndWhile
  \State $i \gets i-1$
  \EndWhile
\State $k \gets k-1$
\EndWhile
\State $\widehat{C}(1,1)=0,\widehat{C}(m,m)=0$
\State \textbf{Output} $\widehat{C}$
\end{algorithmic}
\end{breakablealgorithm}

For the Kaczmarz-Tanabe's iteration and the symmetric Kaczmarz-Tanabe's iteration, the matrices $C$ and $\bar{C}$ are invariant in the subsequent iterations which is beneficial for computation, e.g., in medical imaging equipments, one can calculate and store the matrices $C$ and $\bar{C}$ or related matrices in the imaging device in advance. $C$ and $\bar{C}$ can be calculated by block mode or parallel block mode, which will greatly reduce the cost to compute them. Blocking technology can be made on the linear system which has been discussed in some articles (please refer to \cite{Needell2014paved,Ma2015Convergence,Elfving1980,Censor2002,Needell2015Randomized} for more details).

\section{Numerical tests}
\label{section.numerical.tests}
We will test the convergence rates of the Kaczmarz-Tanabe type methods and compare them with the SIRT and CGMN methods with two examples.
Let $\{y_k,k>0\}$ be the iterative sequences of these methods, and we mainly consider three kinds of iterative errors, i.e., $\|y_k-x^*\|_2$, $\|y_k-x^\dagger\|_2$, and $\|y_k-x^\dagger-P_{N(A)}x_0\|_2$.

For the Kaczmarz-Tanabe methods, the pre-calculation cost of $C$ is $O(m^4)$, the calculation cost of the Kaczmarz-Tanabe's iteration is $O(m^2n)$, and so is the symmetric Kaczmarz-Tanabe method. For the Kaczmarz method, the calculation cost of the Kaczmarz's iteration repeated $m$ times is $O(mn)$.

In addition, $A^TC^TM$ and $A^TC^TMA$ in the Kaczmarz-Tanabe's iteration can also be pre-calculated. Regardless of the pre-calculation cost, the calculation amount of pure Kaczmarz-Tanabe's iteration is only $O(n^2)$. In the sense of pre-calculation, the Kaczmarz-Tanabe type methods are particularly suitable for the over-determined systems with the same projective matrix and different measurement vectors $b$.

\subsection{Tanabe's problem}
Consider the following linear system with equations
\begin{align}\label{Tanabe.problem.1}
  \left (
     \begin{array}{rrrr}
         1.0 &  3.0 &  2.0 & -1.0\\
         1.0 &  2.0 & -1.0 & -2.0\\
         1.0 & -1.0 &  2.0 &  3.0\\
         2.0 &  1.0 &  1.0 &  1.0\\
         5.0 &  5.0 &  4.0 &  1.0\\
         4.0 & -1.0 &  5.0 &  7.0
      \end{array}
  \right  )x=
  \left (
     \begin{array}{r}
        5.0\\
        0.0\\
        5.0\\
        5.0\\
       15.0\\
       15.0
     \end{array}
  \right ).
\end{align}
Linear system \eqref{Tanabe.problem.1} is consistent and over-determined. The general solution is
\begin{align}\label{general.soluton.model1}
  (x_1,x_2,x_3,x_4)^T=k(-2/3,1,-2/3,1)^T+(5/3,0,5/3,0)^T,
\end{align}
where $k\in \mathbb{C}$ is any constant and $\mathbb{C}$ is the complex field. In numerical experiments, $x^*=(1,1,1,1)^T$ is taken as the test solution. We compare the convergence rates of the Kaczmarz-Tanabe and symmetric Kaczmarz-Tanabe methods on the one hand, and compare those of the Kaczmarz-Tanabe type methods and SIRT methods on the other hand.

\begin{figure}[!hbtp]
   \centering
   \begin{minipage}[ht]{0.9\linewidth}
     \centering
     \subfigure[]{
     \includegraphics[width=0.35\linewidth]{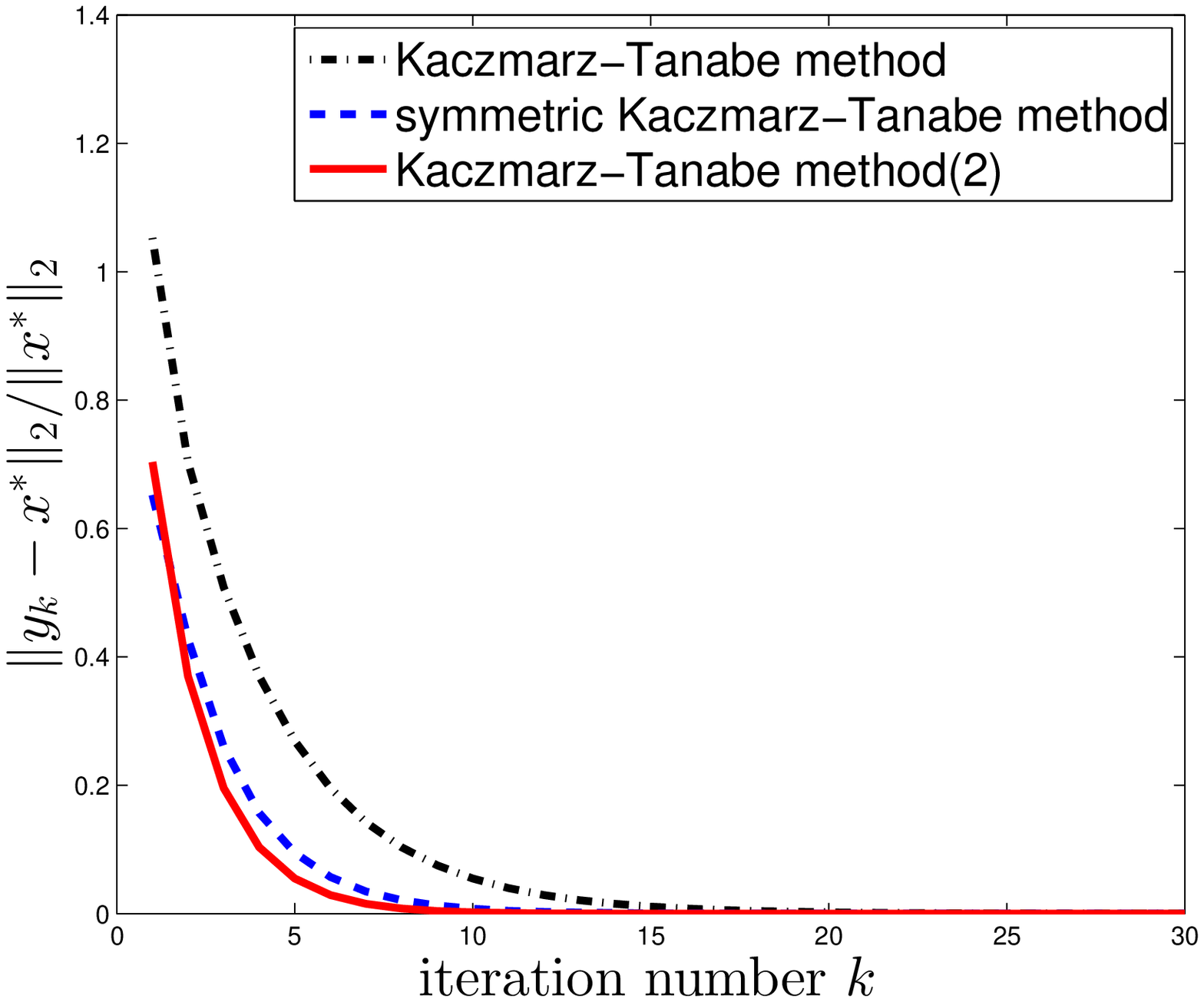}}
     \subfigure[]{
     \includegraphics[width=0.35\linewidth]{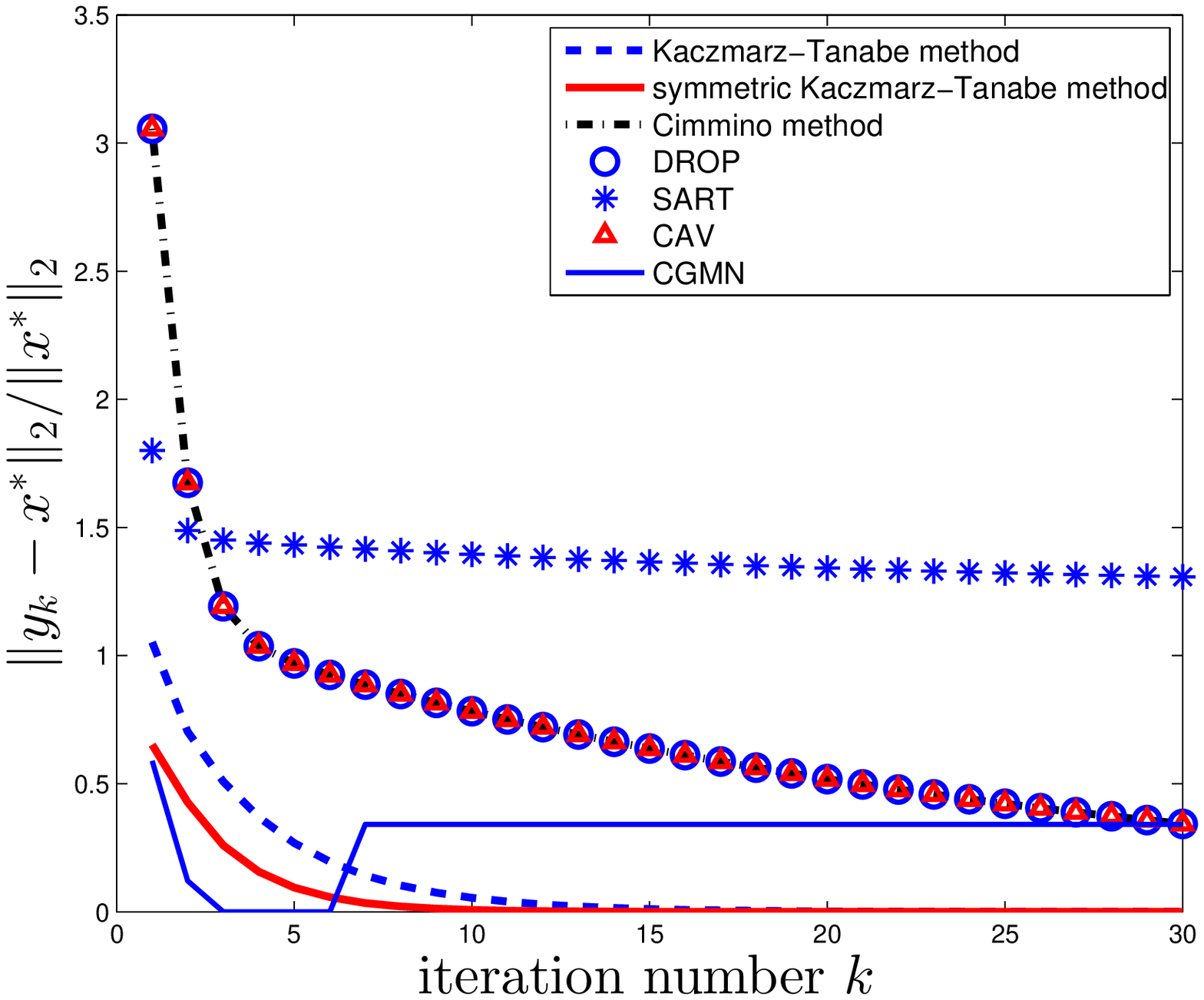}}
  \end{minipage}
%
   \centering
   \begin{minipage}[ht]{0.9\linewidth}
     \centering
     \subfigure[]{
     \includegraphics[width=0.35\linewidth]{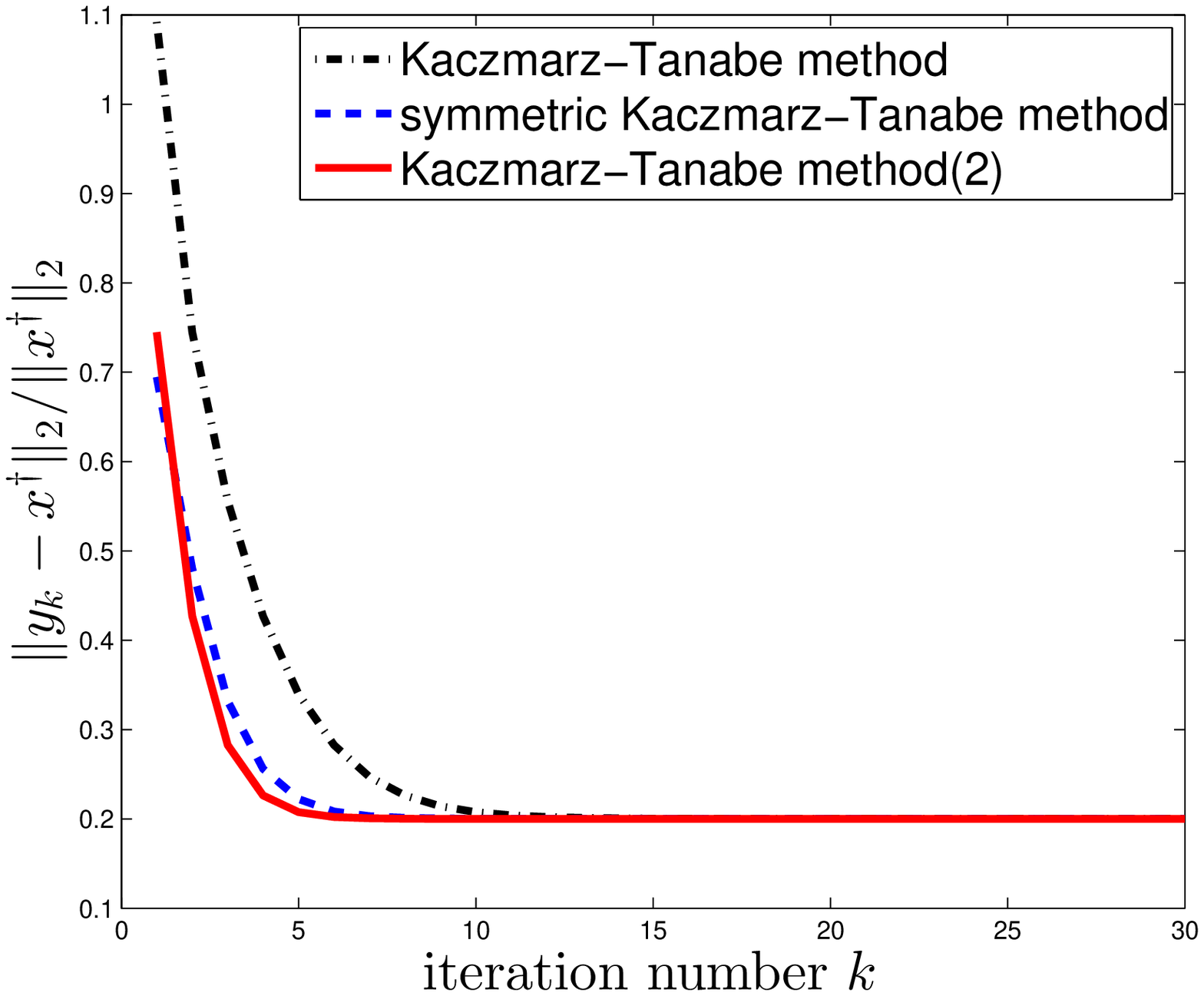}}
     \subfigure[]{
     \includegraphics[width=0.35\linewidth]{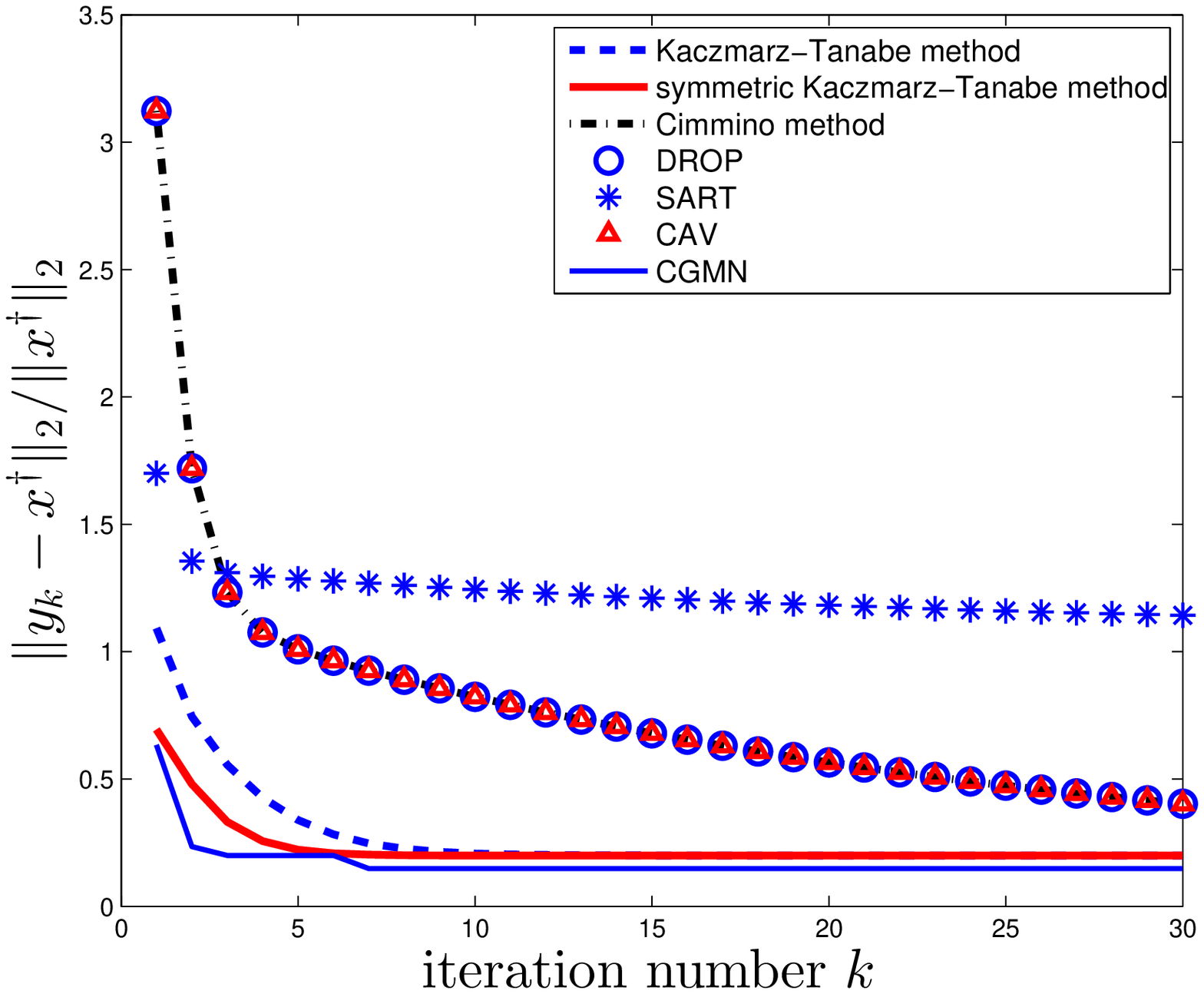}}
  \end{minipage}
%
  \centering
  \begin{minipage}[ht]{0.9\linewidth}
     \centering
     \subfigure[]{
     \includegraphics[width=0.35\linewidth]{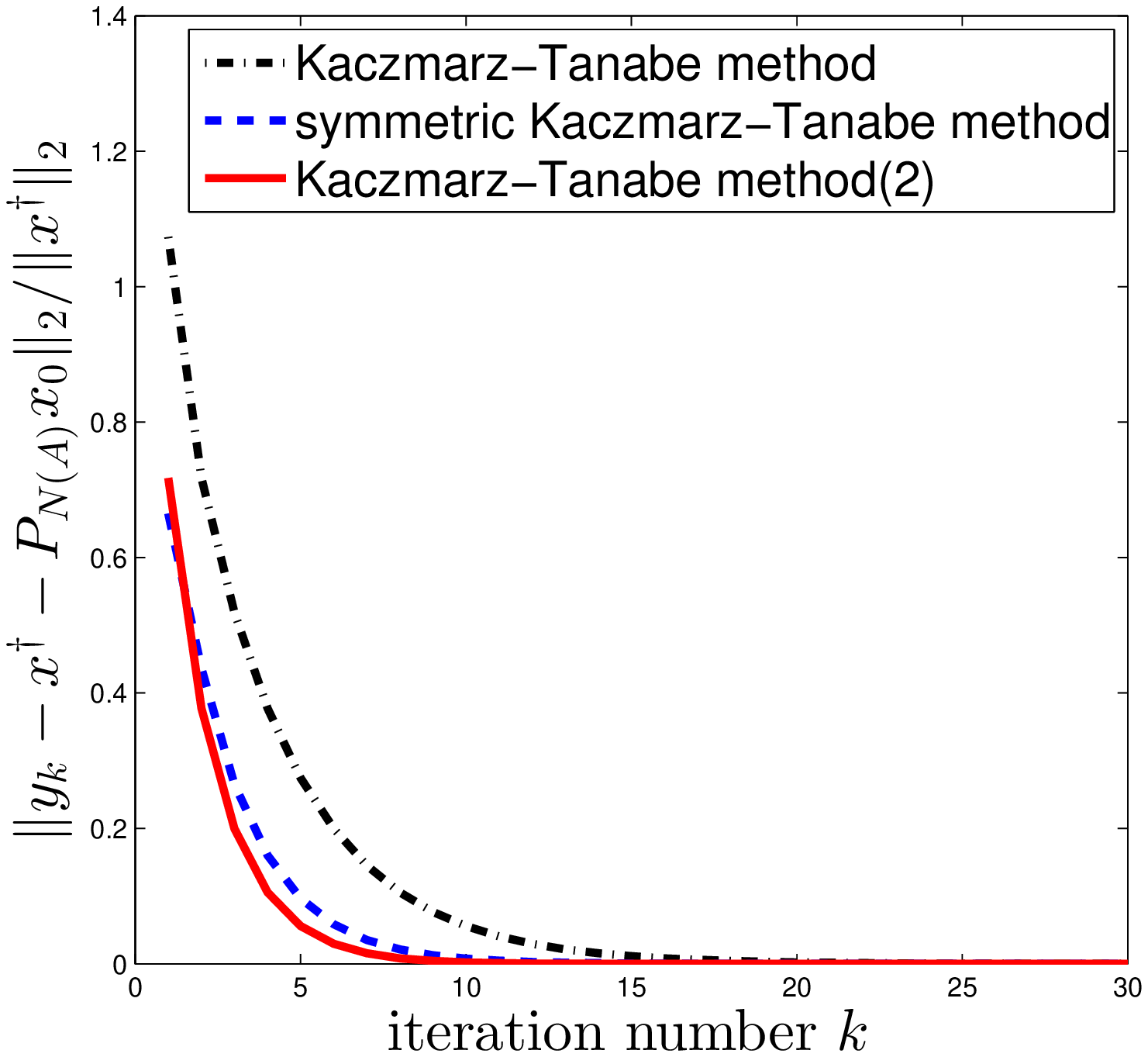}}
     \subfigure[]{
     \includegraphics[width=0.35\linewidth]{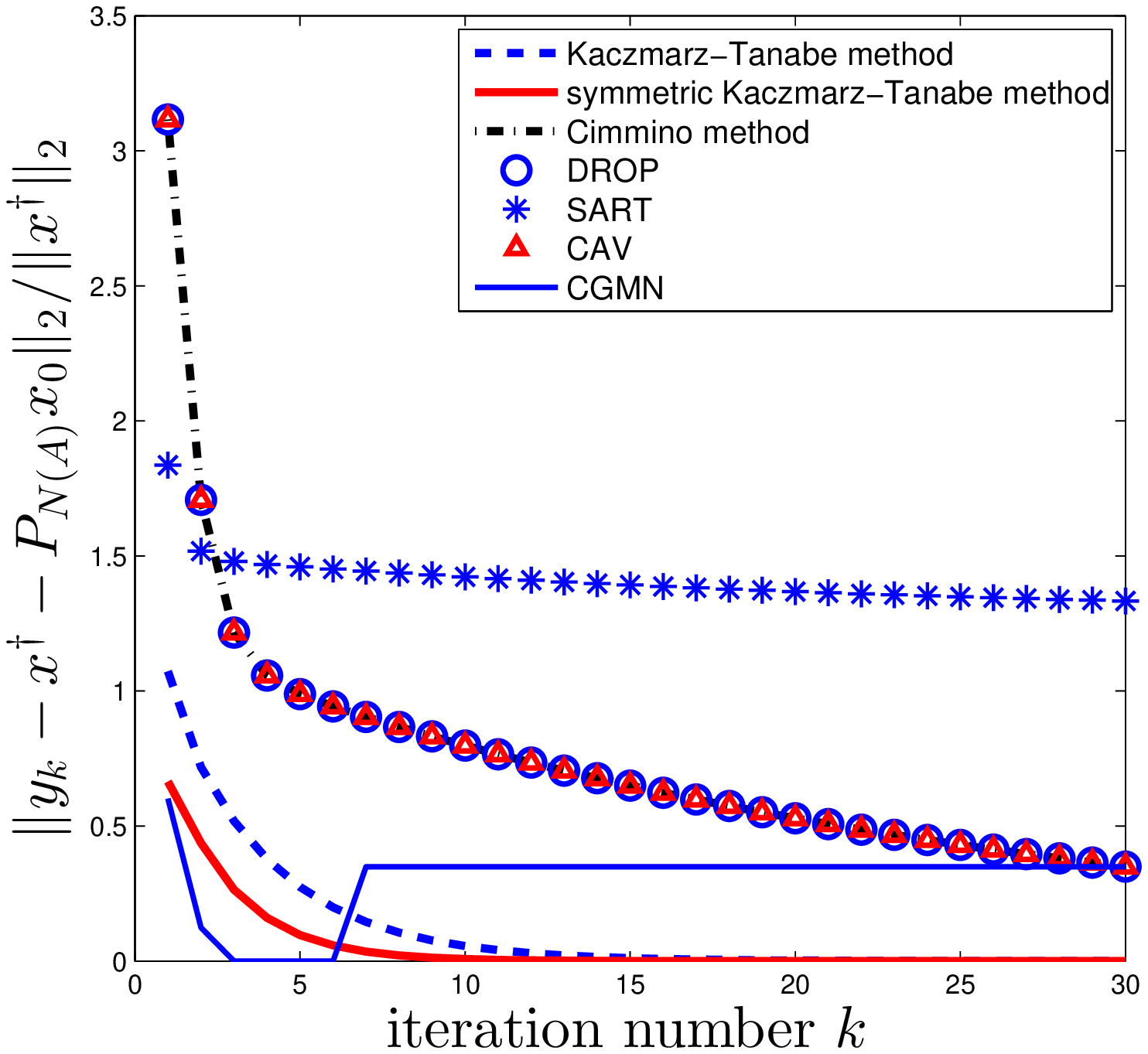}}

  \caption{The comparisons of $\|y_k-x^*\|_2$, $\|y_k-x^\dagger\|_2$ and $\|y_k-x^\dagger-P_{N(A)}x_0\|_2$ when $x_0=(7,6,10,6)^T$, where (a),(c) and (e) are comparisons among the Kaczmarz-Tanabe method, the symmetric Kaczmarz-Tanabe method and two-step Kaczmarz-Tanabe method for solving Tanabe's problem, and (b), (d) and (f) are comparisons among the Kaczmarz-Tanabe method, the symmetric Kaczmarz-Tanabe method and SIRT type methods for solving Tanabe's problem. (see \eqref{the.SIRT.method} for the iterative schemes).}
\label{figure.comparison.convergence.rate.Tanabe.model.1}
 \end{minipage}
\end{figure}

\begin{figure}[!hbtp]
   \centering
   \begin{minipage}{0.9\linewidth}
     \centering
      \setlength{\abovecaptionskip}{0.cm}
      \setlength{\belowcaptionskip}{0.cm}
       \subfigure[]{
       \includegraphics[width=0.35\linewidth]{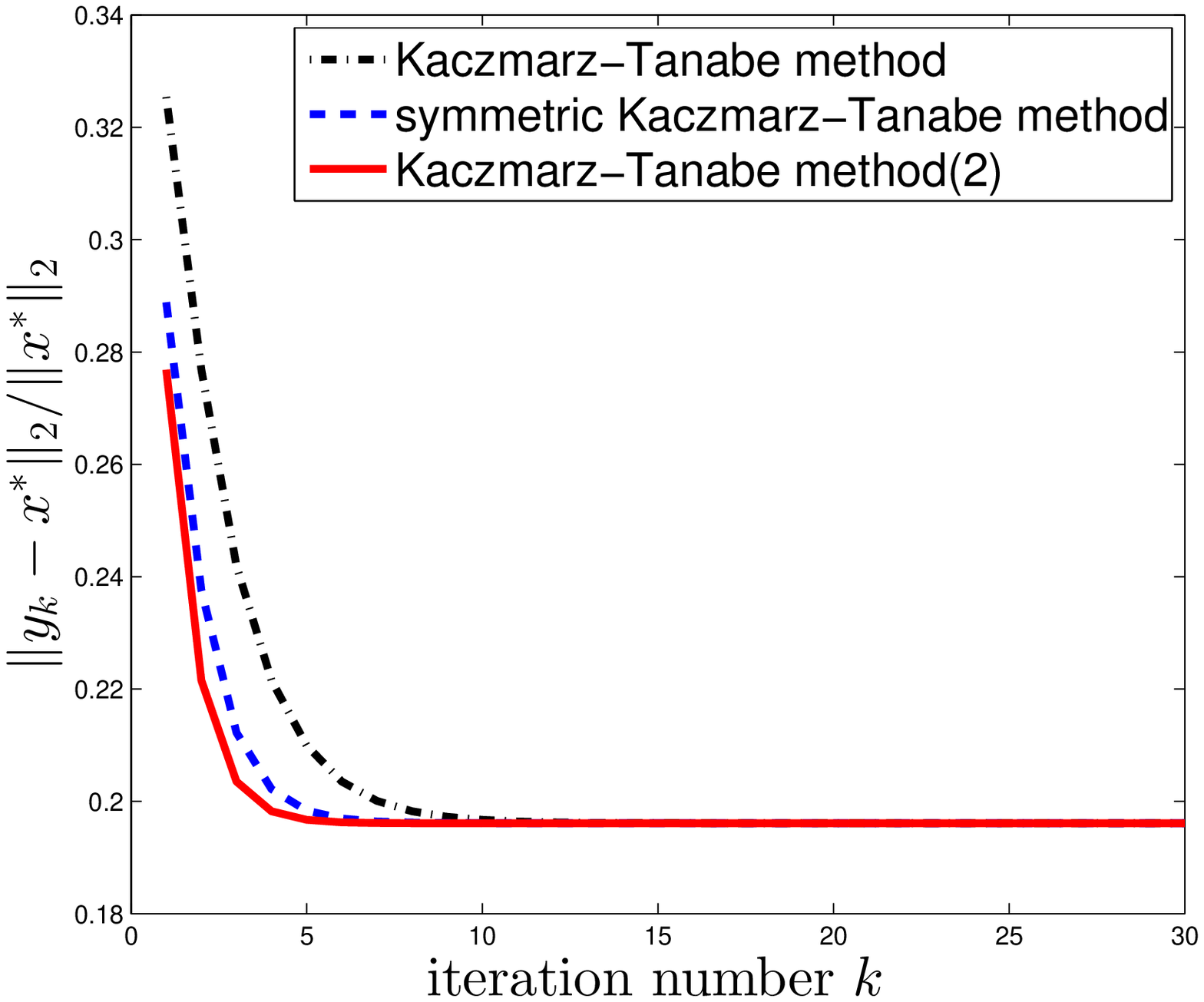}}
       \subfigure[]{
       \includegraphics[width=0.35\linewidth]{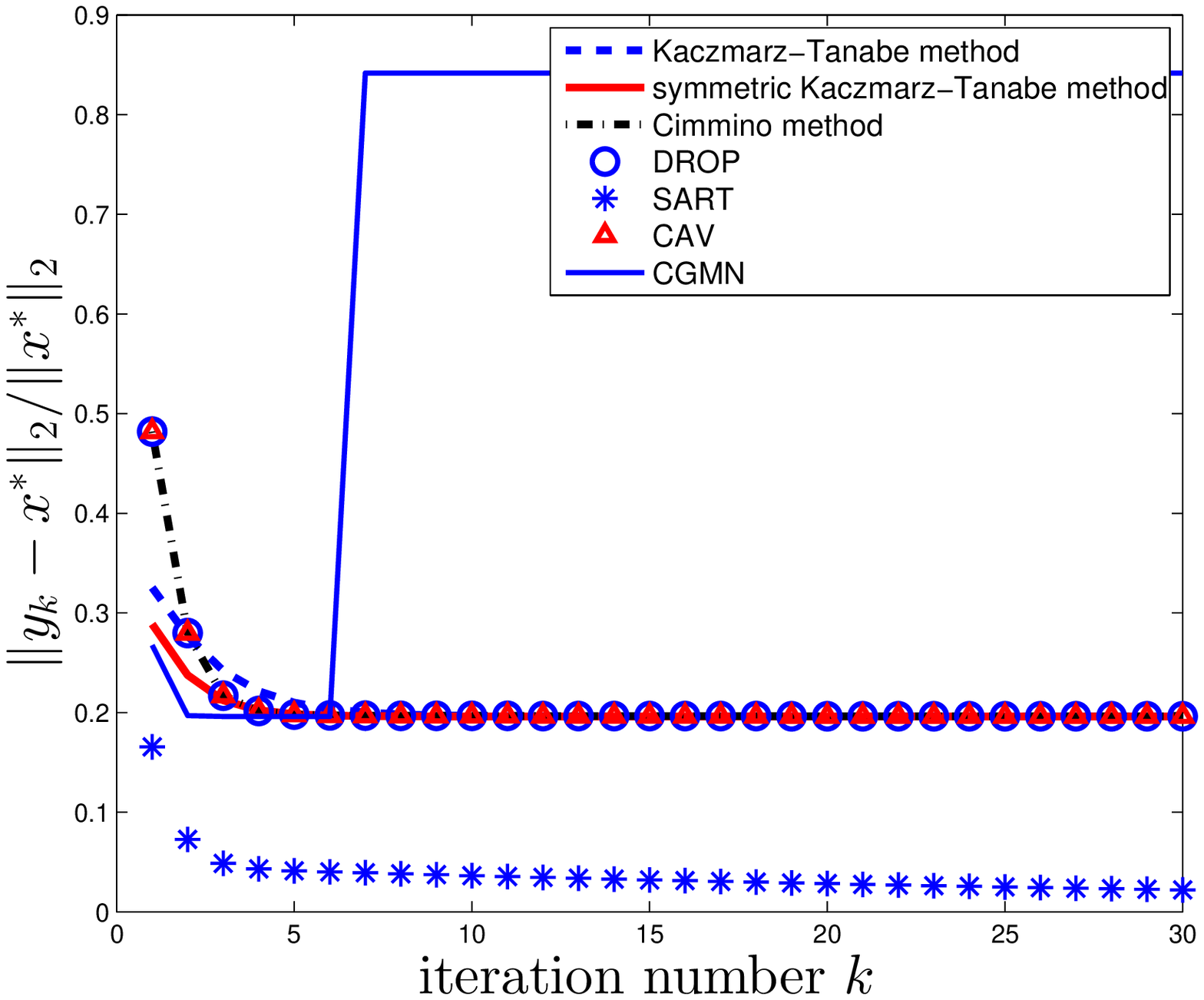}}

      \subfigure[]{
      \includegraphics[width=0.35\linewidth]{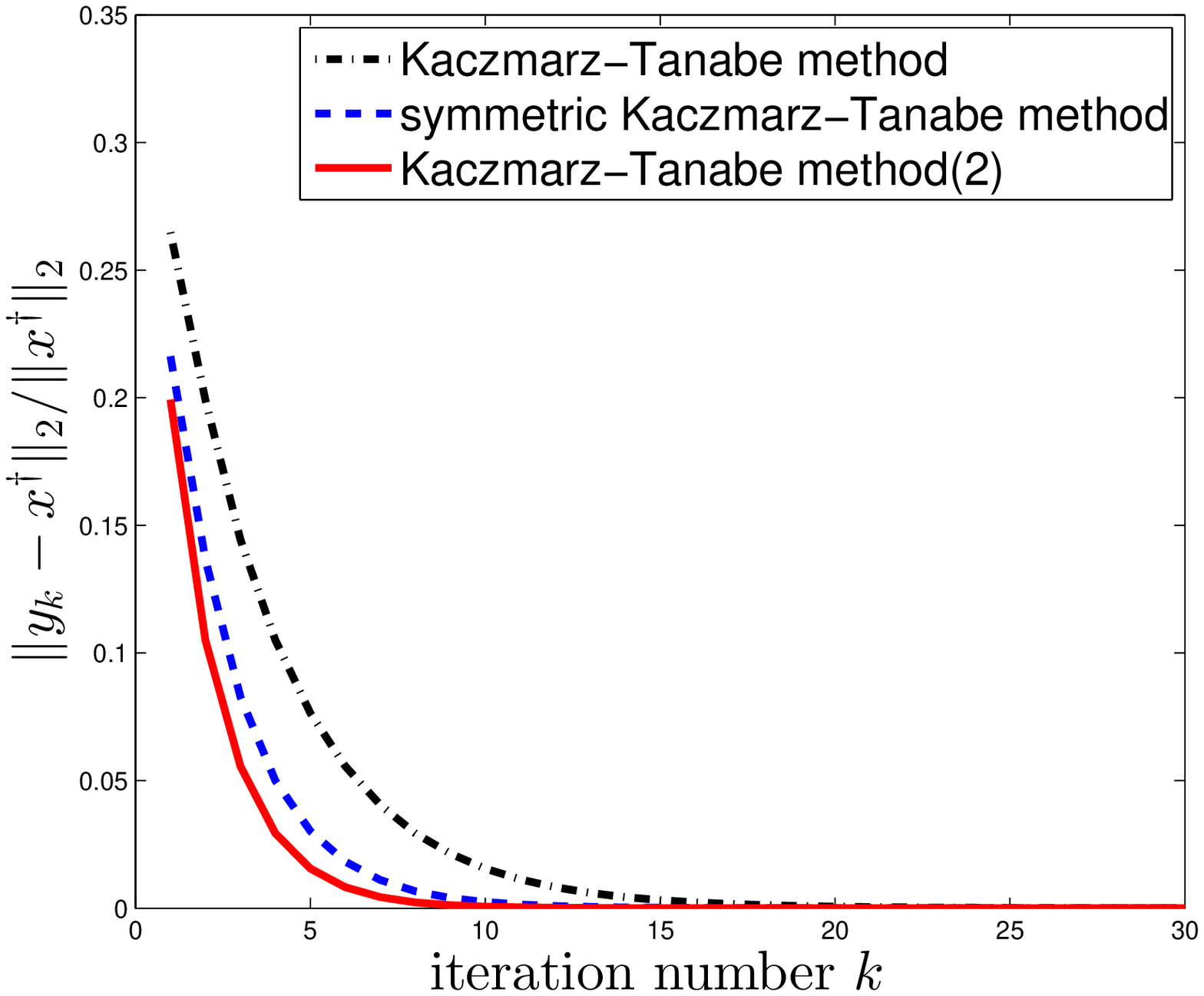}}
      \subfigure[]{
      \includegraphics[width=0.35\linewidth]{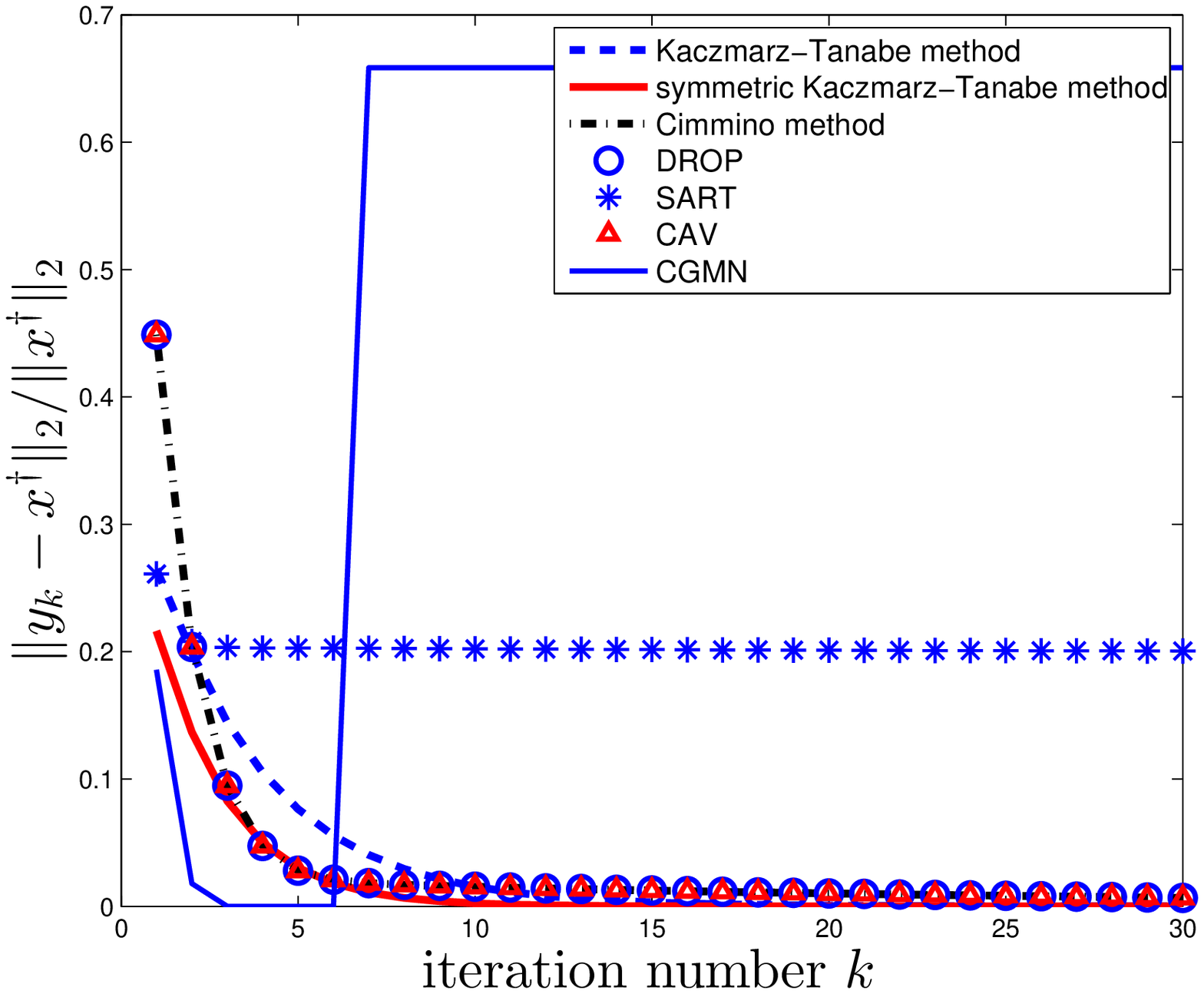}}

      \caption{The comparisons of $\|y_k-x^*\|_2$, $\|y_k-x^\dagger\|_2$ when $x_0=(0,0,0,0)^T$, where (a) and (c) are comparisons among the Kaczmarz-Tanabe method, the symmetric Kaczmarz-Tanabe method and two-step Kaczmarz-Tanabe method for solving Tanabe's problem, and (b) and (d) are comparisons among the Kaczmarz-Tanabe method, the symmetric Kaczmarz-Tanabe method and SIRT type methods for solving Tanabe's problem.}
   \label{figure.comparison.convergence.rate.Tanabe.model.2}
   \end{minipage}
\end{figure}

Numerical results are shown in Figures \ref{figure.comparison.convergence.rate.Tanabe.model.1}$\sim$\ref{figure.comparison.convergence.rate.Tanabe.model.2}, where Figure \ref{figure.comparison.convergence.rate.Tanabe.model.1} shows the error curves of $\|y_k-x^*\|_2$, $\|y_k-x^\dagger\|_2$, and $\|y_k-x^\dagger-P_{N(A)}x_0\|_2$ when $x_0=(7,6,10,6)^T$, and Figure \ref{figure.comparison.convergence.rate.Tanabe.model.2} shows the corresponding results when $x_0=(0,0,0,0)^T$. In Figures \ref{figure.comparison.convergence.rate.Tanabe.model.1}(a)(c)(e) and \ref{figure.comparison.convergence.rate.Tanabe.model.2}(a)(c), we compare the errors of Kaczmarz-Tanabe method, symmetric Kaczmarz-Tanabe method and two-step Kaczmarz-Tanabe method (marked with `Kaczmarz-Tanabe(2)' in these figures).

In Figures \ref{figure.comparison.convergence.rate.Tanabe.model.1}(b)(d)(f) and \ref{figure.comparison.convergence.rate.Tanabe.model.2}(b)(d), we compare the errors of Kaczmarz-Tanabe method, symmetric Kaczmarz-Tanabe method, Cimmino method, DROP method, SART method, CAV method and CGMN method when $x_0=(7,6,10,6)^T$ and $x_0=(0,0,0,0)^T$ respectively. Since the computational work of the Kaczmarz-Tanabe method and the symmetric Kaczmarz-Tanabe method is roughly the same as that of the SIRT methods when $C$ and $\bar{C}$ are determined, therefore we deal with these methods in the same way, that is, comparing one Kaczmarz-Tanabe' iteration with one symmetric Kaczmarz-Tanabe's iteration, as well as other methods.
In Figure \ref{figure.comparison.convergence.rate.Tanabe.model.1}, (a),(b) are the same as (e),(f) respectively, although they look different. Denote
\begin{align*}
  \xi=(-2/3,1,-2/3,1)^T.
\end{align*}
We know from \eqref{general.soluton.model1} that $N(A)=\textbf{span}\{\xi\}$, thus
\begin{align*}
   &P_{N(A)}x_0=P_{N(A)}x^*=\frac{\xi^Tx_0}{\|\xi\|_2^2}\xi=\frac{3}{13}(-2/3,1,-2/3,1)^T,\\
   &x^*=x^\dagger+P_{N(A)}x_0,
\end{align*}
which means that
\begin{align*}
  \|y_k-x^*\|_2=\|y_k-x^\dagger-P_{N(A)}x_0\|_2.
\end{align*}
Therefore, the convergence of the error curves shown in Figure \ref{figure.comparison.convergence.rate.Tanabe.model.1} (a), (b), (e) and (f) are consistent with the theoretical results, and this is also why the curves in Figure \ref{figure.comparison.convergence.rate.Tanabe.model.1} (c) and (d) do not tend to the $x$-axis.

In addition, Figure \ref{figure.comparison.convergence.rate.Tanabe.model.1} (a), (c) and (e) also show that one symmetric Kaczmarz-Tanabe's iteration is better than one Kaczmarz-Tanabe's iteration, and  slightly worse than the two-step Kaczmarz-Tanabe's iteration. Meanwhile, Figure \ref{figure.comparison.convergence.rate.Tanabe.model.1} (b), (d) and (f) show that the convergence speed of the Kaczmarz-Tanabe and symmetric Kaczmarz-Tanabe methods is faster than the SIRT methods.

Figure \ref{figure.comparison.convergence.rate.Tanabe.model.2} shows the efficiency of these methods when $x_0=(0,0,0,0)^T$. Figure \ref{figure.comparison.convergence.rate.Tanabe.model.2} (a) is slightly different from Figure \ref{figure.comparison.convergence.rate.Tanabe.model.1} (a), and Figure \ref{figure.comparison.convergence.rate.Tanabe.model.2} (c) is consistent with Figure \ref{figure.comparison.convergence.rate.Tanabe.model.1} (e). It seems from Figure \ref{figure.comparison.convergence.rate.Tanabe.model.2} (b) that the SART method is better than the others. The reason is that the SART's iteration converges to $x^*$ rather than $x^\dagger$ when $x_0=(0,0,0,0)^T$, which can be seen from Figure \ref{figure.comparison.convergence.rate.Tanabe.model.2} (d).

We also note that the CGMN method is sensitive to iteration step, and converges quickly at the beginning, and then the results become worse. Suppose the linear system to be solved by CGMN method is $Bx=c$, this phenomenon may be related to the positive semi-definiteness of $B$. In other words, the descending direction $d$ of the conjugate gradient (CG) method becomes an eigenvector of $0$ eigenvalue of $B$ or $Bd\approx 0$.

\subsection{Headphantom problem}
In computerized tomography, the distribution of some physical parameter(such as absorption intensities) at the cross-section of the object need to be reconstructed from the projection data such as medical diagnosis--the distribution of the absorption intensities of tissue slice need to be reconstructed from {X}-ray data. The computerized tomography system attributes to a linear system $Ax=b$, where $A$ is a projected system, $b$ is scanning data, and $x$ is unknown intensity image of an object. In the general case, the system is overdetermined.

The linear system is generated from the subroutine `parallel' in AIR Tool II package \cite{Hansen2018AIRtools}, and there are 36 projective angles at equal intervals in $[0,2\pi]$ and 75 equi-spaced parallel rays per angle. The headphantom is discretized into $50\times 50$ pixels. and the dimension of $A$ is $2700\times 2500$.

\begin{figure}[!hbtp]
   \centering
   \begin{minipage}[ht]{0.9\linewidth}
      \centering
      \setlength{\abovecaptionskip}{0.cm}
      \setlength{\belowcaptionskip}{0.cm}
      \subfigure[]{
      \includegraphics[width=0.35\linewidth]{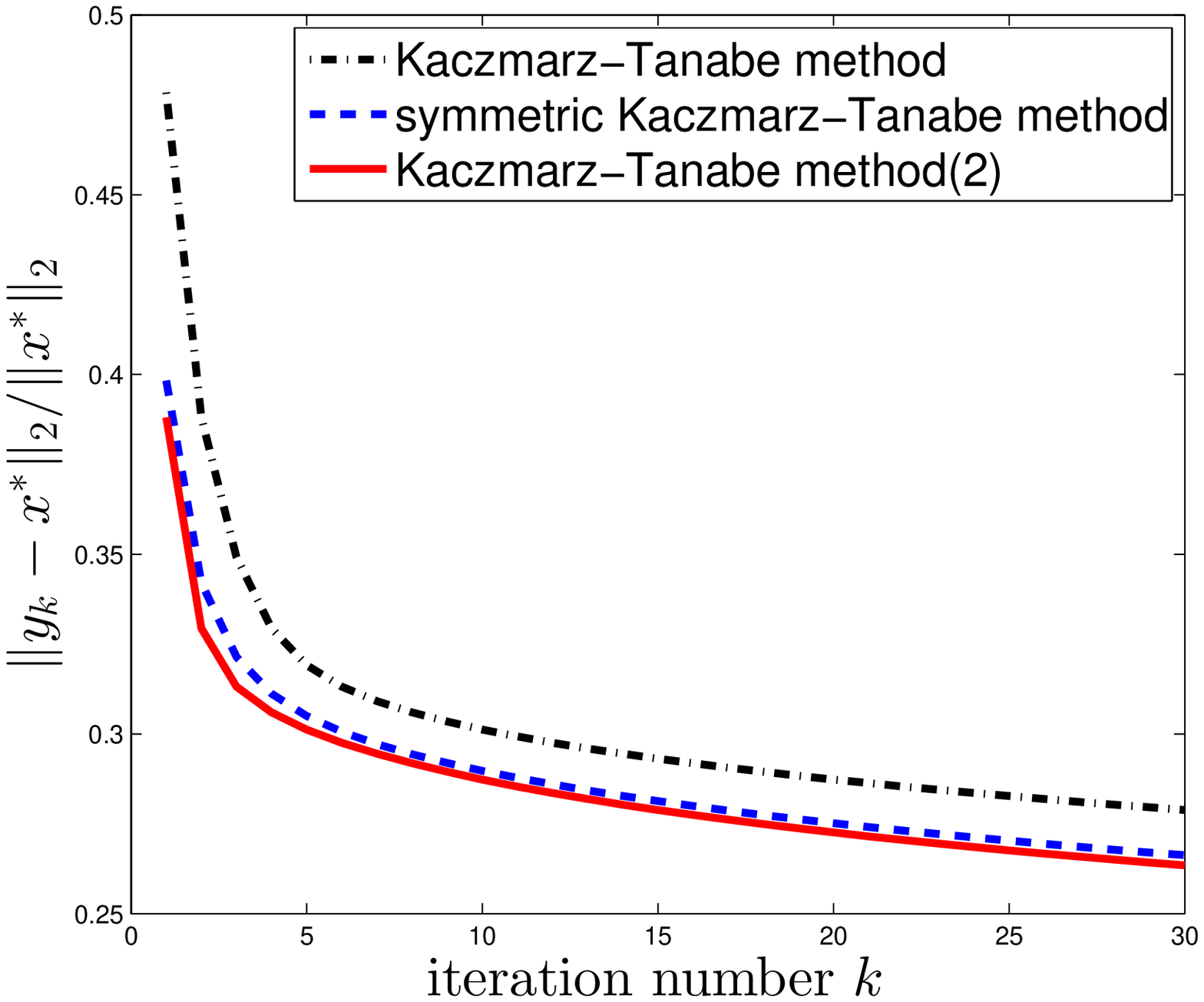}}
      \subfigure[]{
      \includegraphics[width=0.35\linewidth]{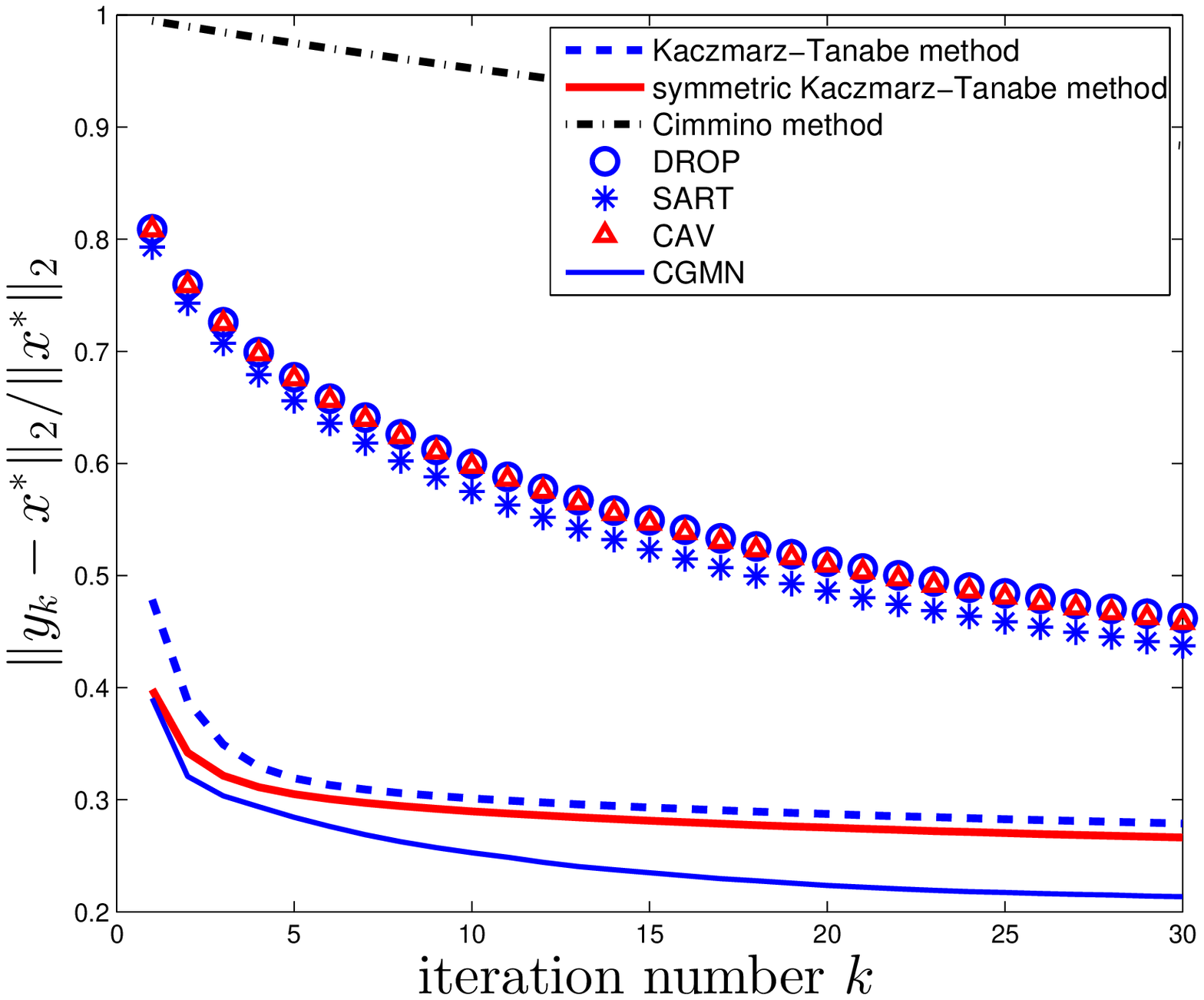}}

      \subfigure[]{
      \includegraphics[width=0.35\linewidth]{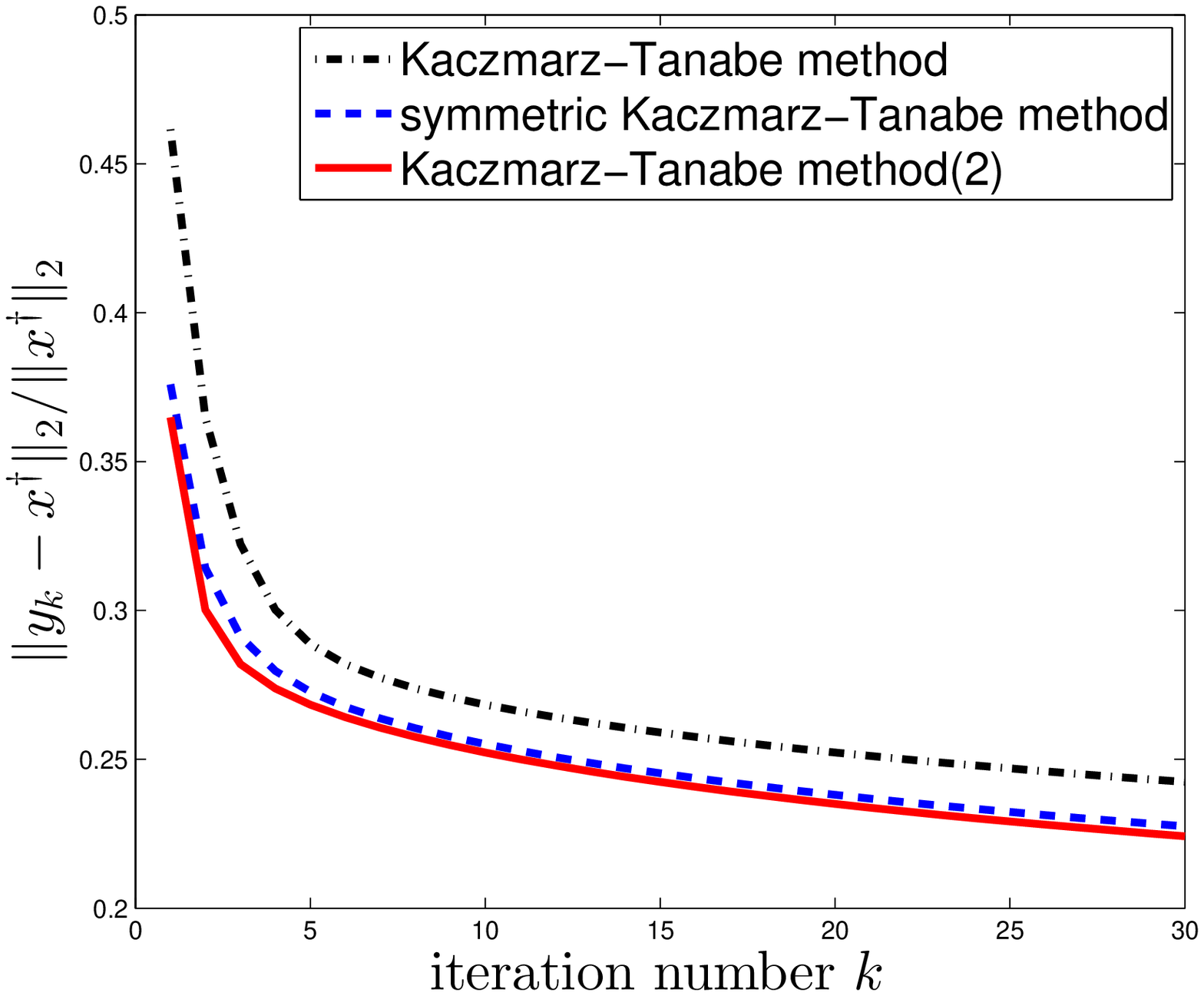}}
      \subfigure[]{
      \includegraphics[width=0.35\linewidth]{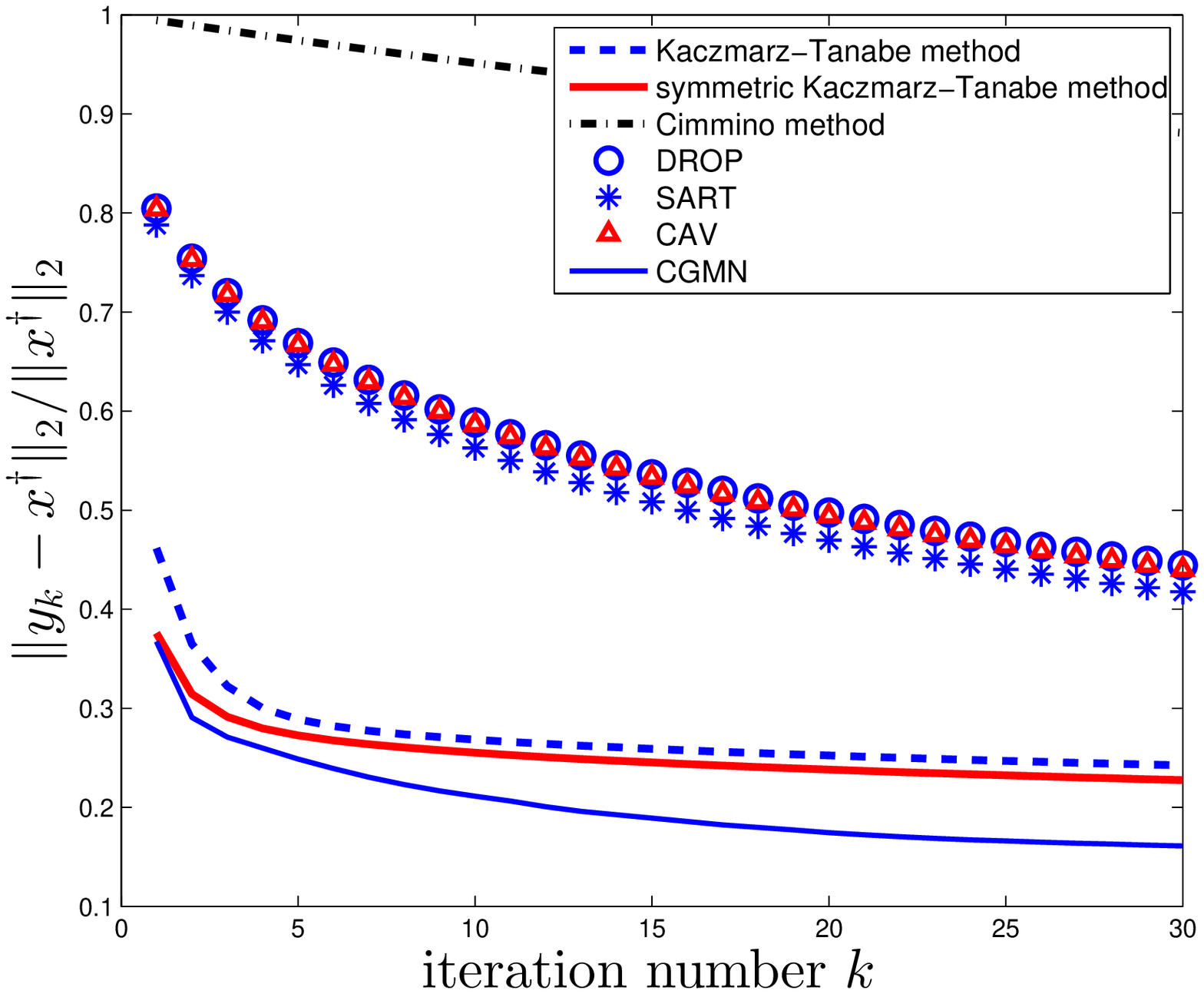}}

      \caption{The comparisons of $\|y_k-x^*\|_2$, $\|y_k-x^\dagger\|_2$, where (a) and (c) are comparisons among the Kaczmarz-Tanabe method, the symmetric Kaczmarz-Tanabe method and two-step Kaczmarz-Tanabe method for solving Headphantom problem, and (b) and (d) are comparisons among the Kaczmarz-Tanabe method, the symmetric Kaczmarz-Tanabe method and SIRT type methods for solving Headphantom problem.}
      \label{figure.comparison.convergence.rate.headphantom.problem}
   \end{minipage}
\end{figure}

The initial value is taken as $x_0=\bf{0}\in \mathbb{R}^{2500}$, and numerical results are shown in Figure \ref{figure.comparison.convergence.rate.headphantom.problem}, where (a) and (c) are results of the Kaczmarz-Tanabe method, symmetric Kaczmarz-Tanabe method and two-step Kaczmarz-Tanabe method for solving the Headphantom problem, (b) and (d) are results of the Kaczmarz-Tanabe method, symmetric Kaczmarz-Tanabe method, SIRT type methods and CGMN method for solving the problem. For this problem, the CGMN method seems to be better than the other methods and the phenomenon in the Tanabe's problem does not appear.

\begin{figure}[!hbt]
  \centering
  \begin{minipage}[ht]{.8\linewidth}
      \centering
      \setlength{\abovecaptionskip}{0.cm}
      \setlength{\belowcaptionskip}{0.cm}
      \subfigure[Kaczmarz-Tanabe method]{
      \includegraphics[width=0.28\linewidth]{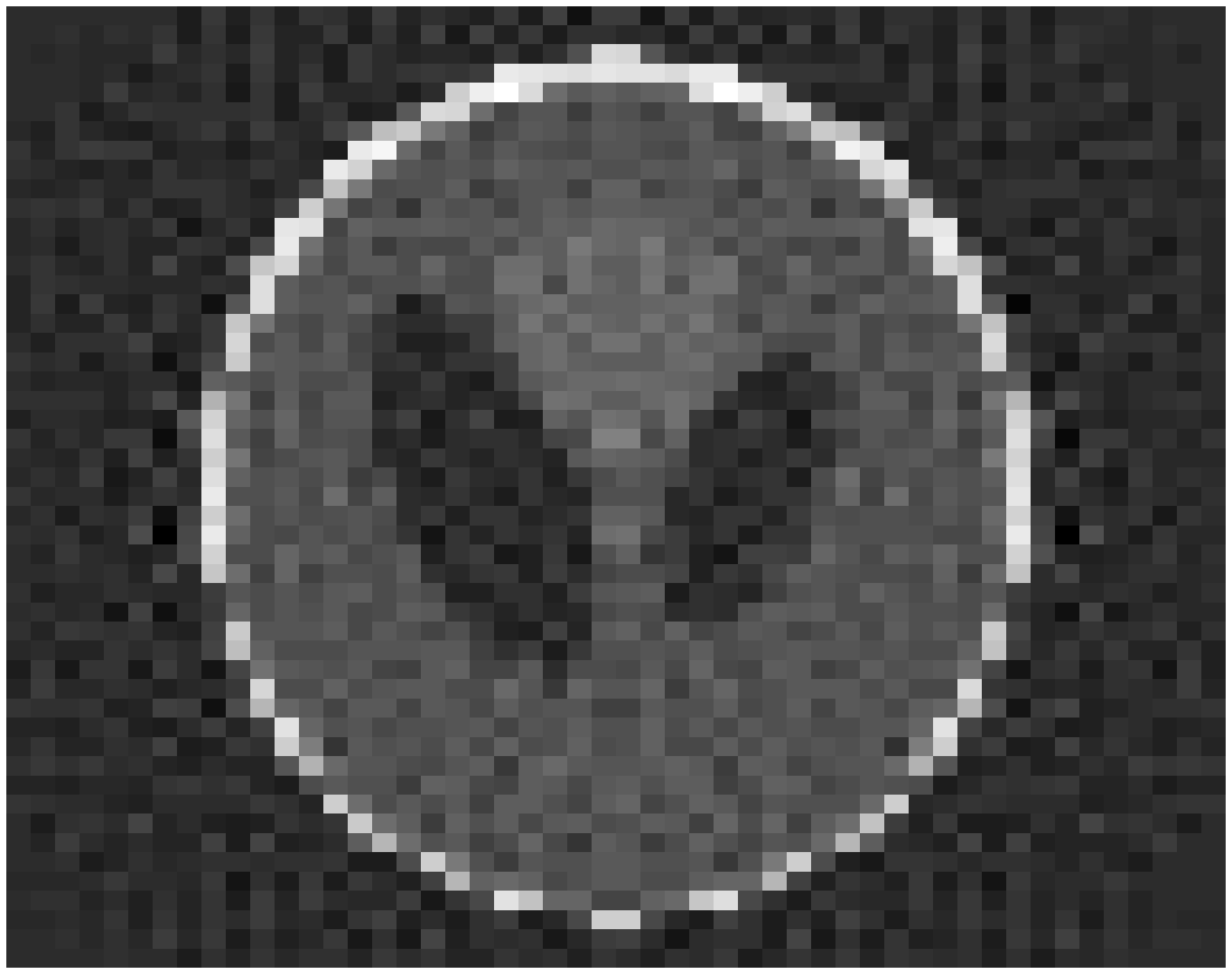}}
      \subfigure[Symmetric Kaczmarz-Tanabe method] {
      \includegraphics[width=0.28\linewidth]{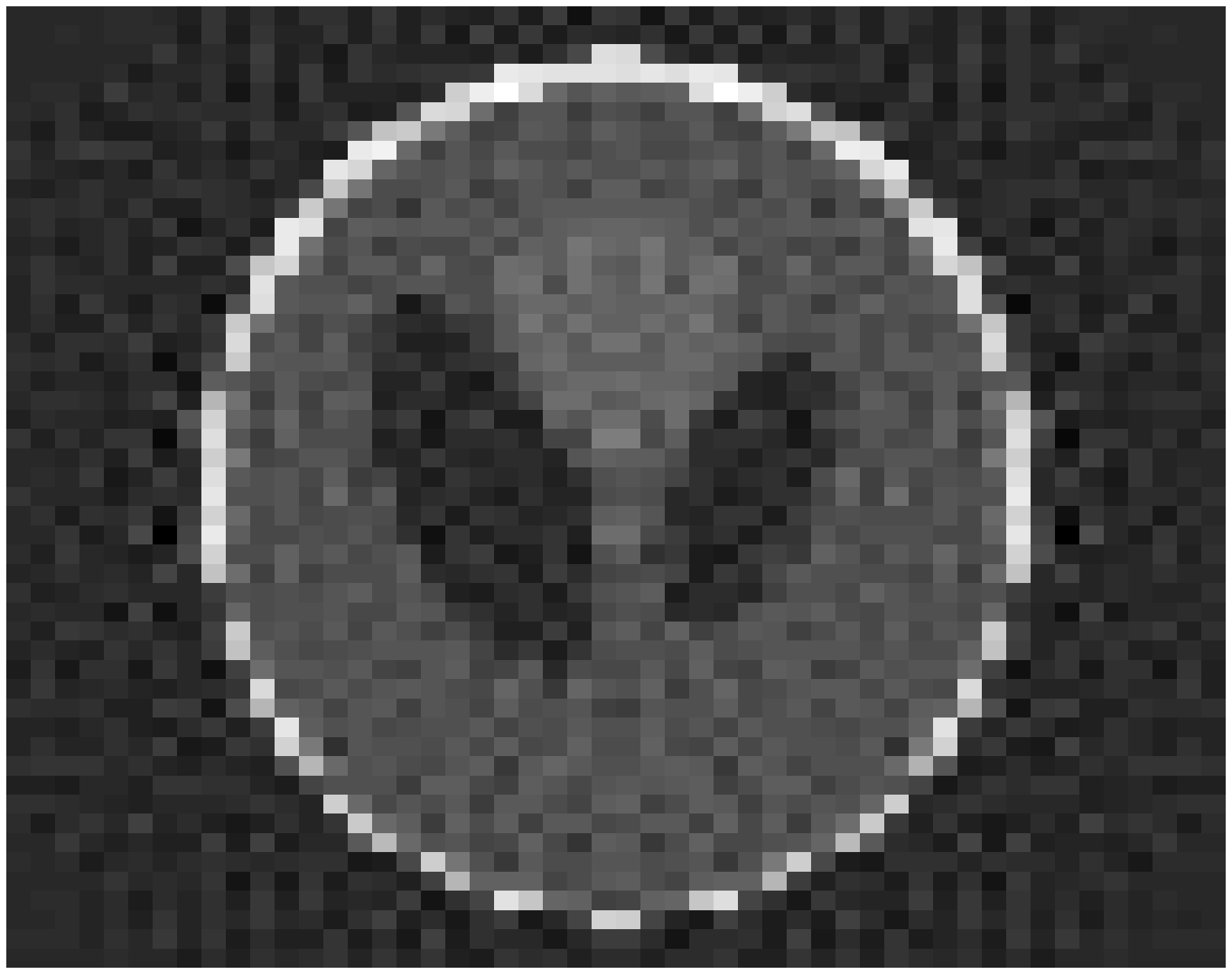}}
      \subfigure[DROP method]{
      \includegraphics[width=0.28\linewidth]{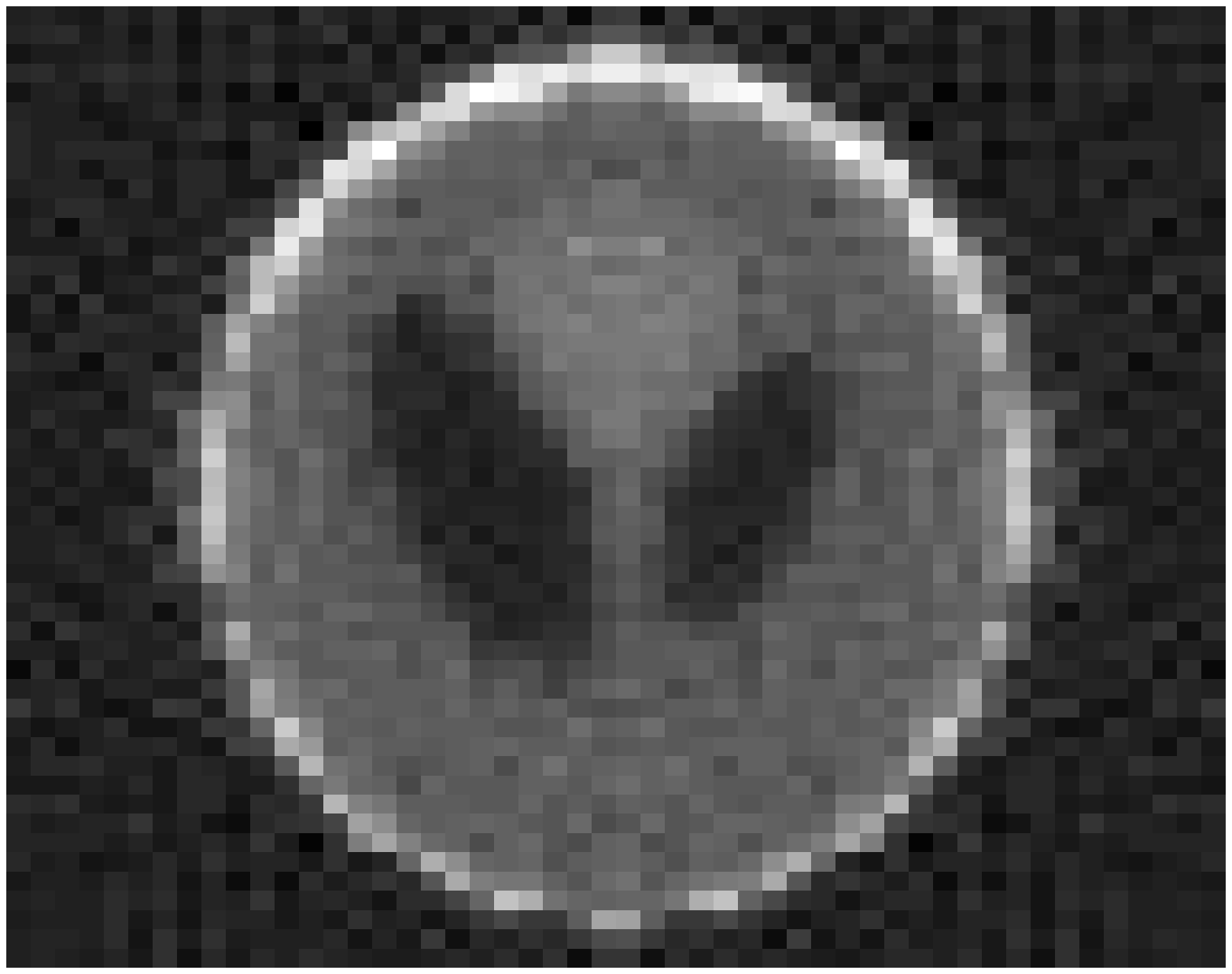}}

      \subfigure[SART method]{
      \includegraphics[width=0.28\linewidth]{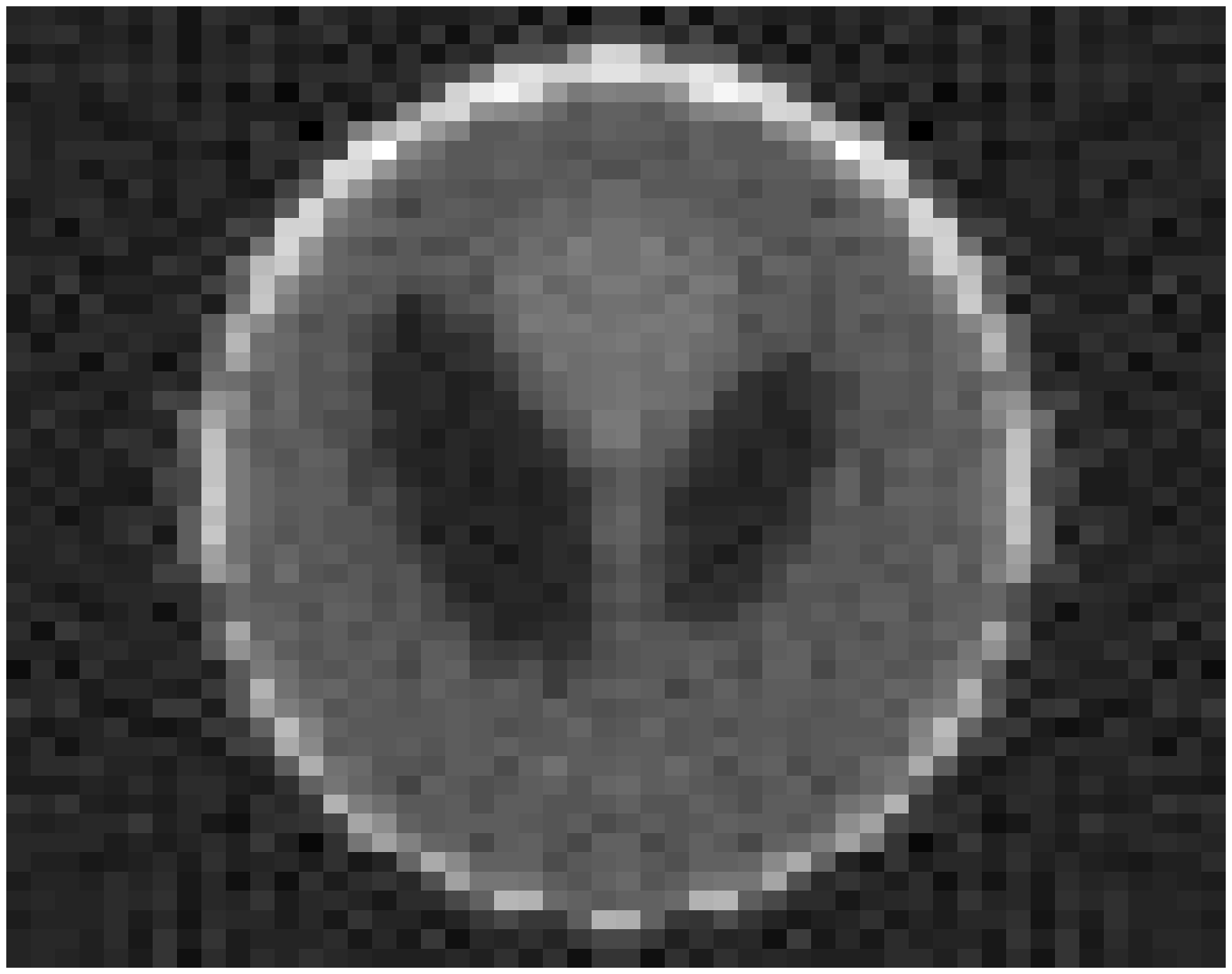}}
      \subfigure[CAV method]{
      \includegraphics[width=0.28\linewidth]{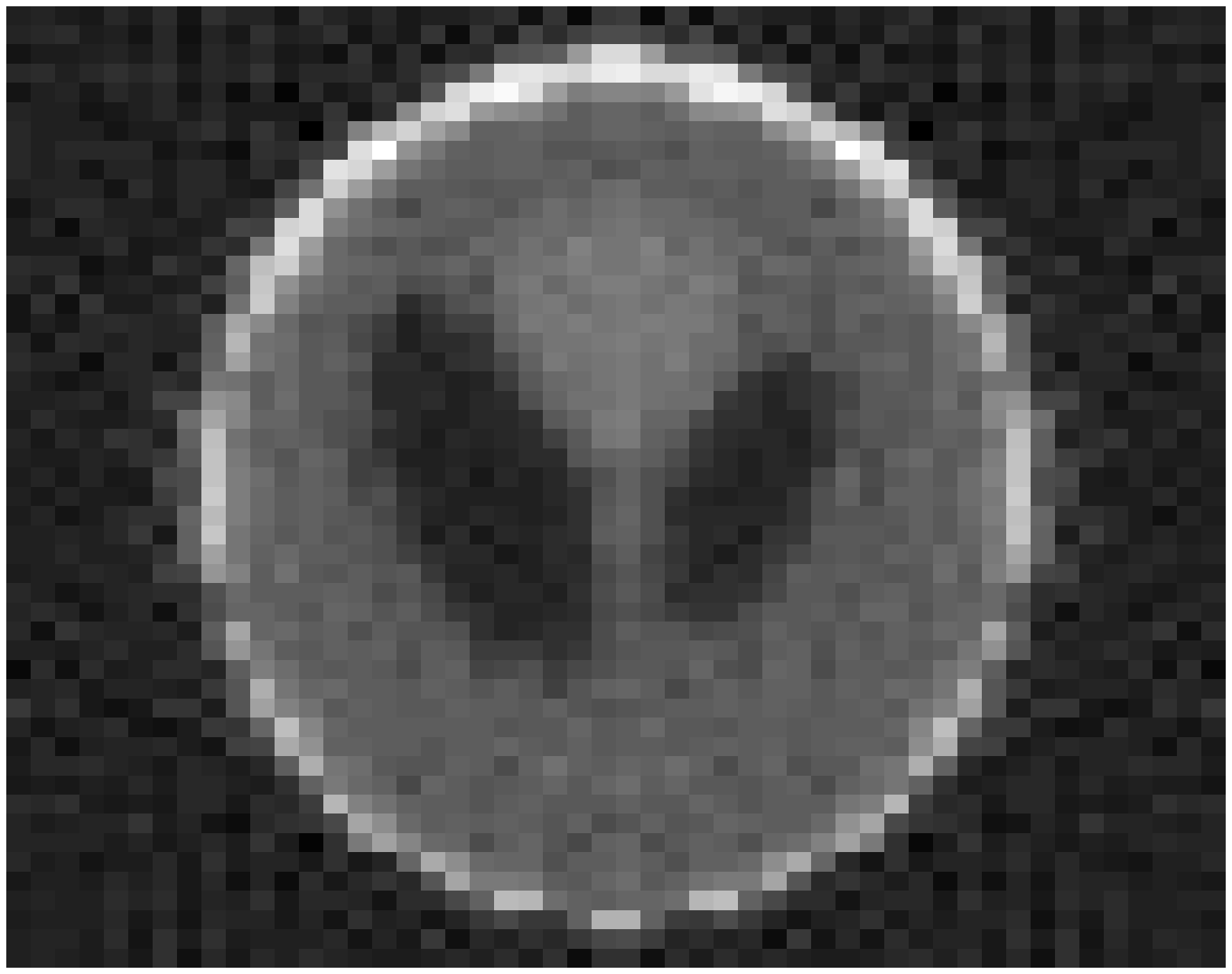}}
      \subfigure[Cimmino method]{
      \includegraphics[width=0.28\linewidth]{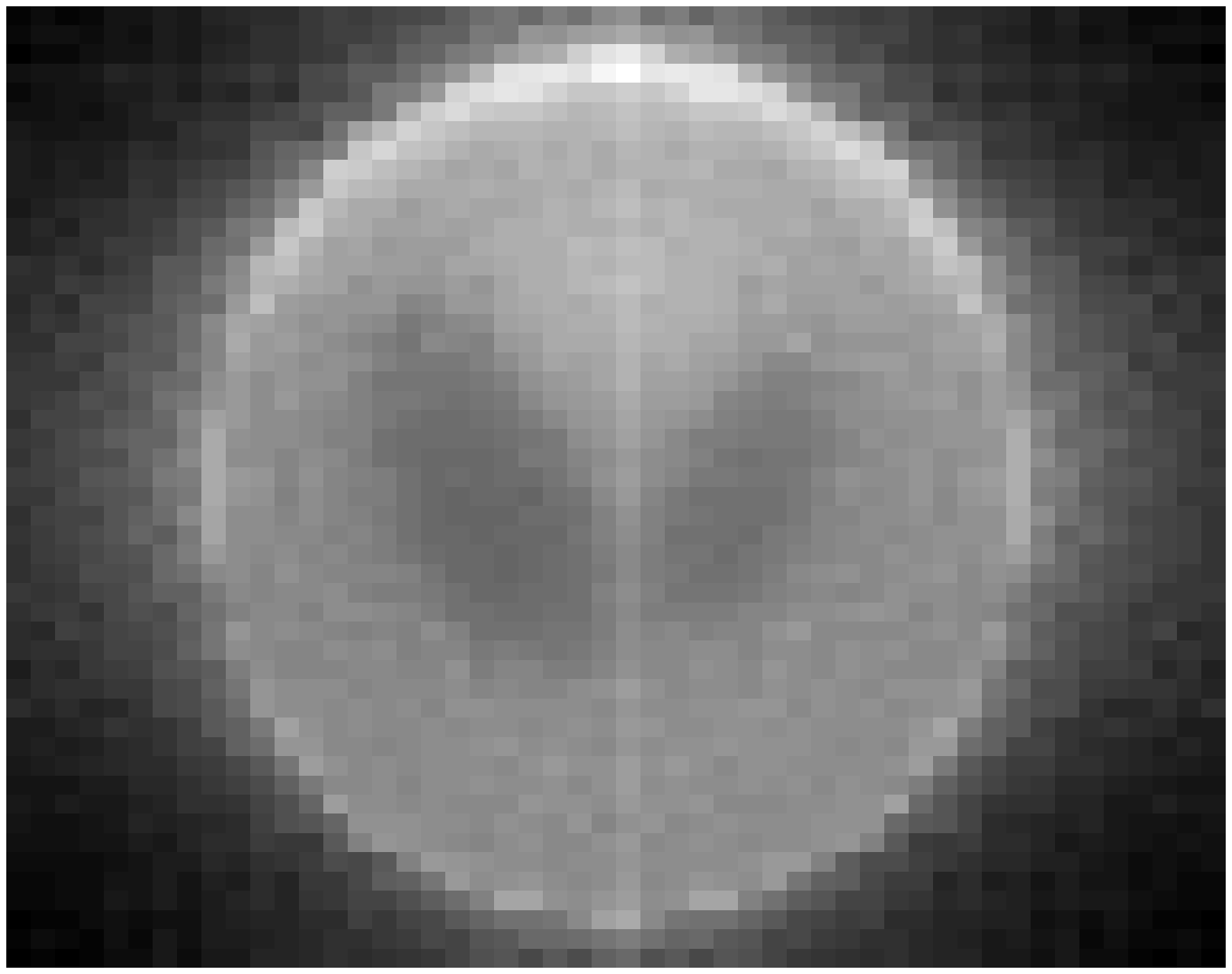}}

      \subfigure[CGMN method]{
      \includegraphics[width=0.28\linewidth]{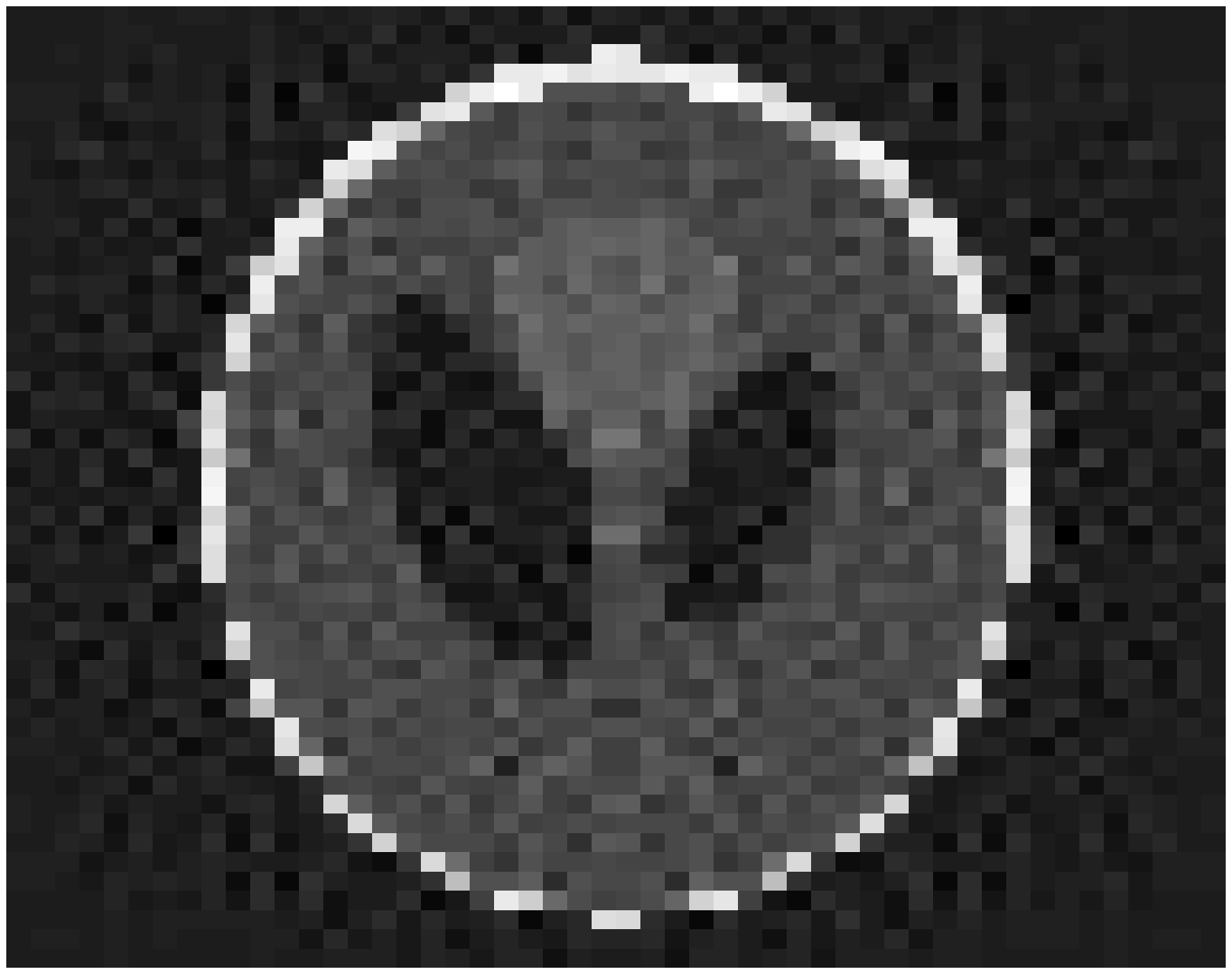}}
      \caption{\footnotesize Numerical images of the Kaczmarz-Tanabe method, symmetric Kaczmarz-Tanabe method and the SIRT type methods for solving Headphantom problem, including DROP, SART, CAV, Cimmino and CGMN methods.}
  \label{images.headphantom.problem}
  \end{minipage}
\end{figure}

From Figure \ref{figure.comparison.convergence.rate.headphantom.problem}, we can see that the Kaczmarz-Tanabe and symmetric Kaczmarz-Tanabe methods are significantly better than the SIRT methods, and slight worse than CGMN method. Numerical images of these methods are shown in Figure \ref{images.headphantom.problem}. From the visual effect, the Kaczmarz-Tanabe method, symmetric Kaczmarz-Tanabe method and CGMN method are close and better than the SIRT type methods.

\section{Conclusion}
\label{section.conclusion}
The Kaczmarz-Tanabe method is the further improvement of the Kaczmarz method. Due to the row to row iterative characteristic of the Kaczmarz method, the Kaczmarz's iteration generally converges slowly and has volatility for perturbed linear systems. The Kaczmarz-Tanabe method overcomes the volatility of Kaczmarz's method and can smoothly approach the `pseudo-inverse' solution when solving the perturbed problem, which lays a foundation for us to further study the minimum norm least-squares solution.

In addition, as a comparison, we also consider the more popular symmetric Kaczmarz-Tanabe method and derive its standard form. We should pay attention to the symmetric Kaczmarz-Tanabe method because one iteration of the symmetric Kaczmarz-Tanabe method can almost obtain the effect of two iterations of the Kaczmarz-Tanabe method. The Kaczmarz-Tanabe's iteration and the symmetric Kaczmarz-Tanabe's iteration have the same iterative formula, if $C$ and $\bar{C}$ are known, then the symmetric Kaczmarz-Tanabe method has obvious advantages over the Kaczmarz-Tanabe method in computational efficiency.

Numerical tests also show that the Kaczmarz-Tanabe type methods, i.e., the Kaczmarz-Tanabe method and the symmetric Kaczmarz-Tanabe method in this paper, are better than the SIRT methods. Although Kaczmarz-Tanabe type methods can not achieve the convergence effect of the CGMN method in some cases, they have advantages in problem applicability, i.e., they converge stably to the minimum norm least-square solution for all compatible linear systems when the initial guess $x_0\in R(A^T)$. In particular, after obtaining $C$ and $\bar{C}$, the Kaczmarz-Tanabe's iteration and the symmetric Kaczmarz-Tanabe's iteration can be implemented as easily as the SIRT methods. In practical applications, such as medical image reconstruction and so on, $C$ and $\bar{C}$ can be calculated in advance and stored in the device, which enables us to implement these iterative methods quickly and get a better solution.

\section{Acknowledgments}
The author thanks the editor and the anonymous reviewers for their constructive comments that helped improve the overall quality and readability of the paper. The research is partially supported by National Science Foundation of China (NSFC) No. 12271401.

\section*{References}
\bibliography{Reference(simp)}

\end{document}